\documentclass [12pt]{article}
\usepackage{latexsym,amssymb,amsfonts,amsmath,amsthm,graphicx}
\pdfoutput=1
\usepackage{graphicx}
\usepackage{hyperref}

\newtheorem{lemma}{Lemma}[section]
\newtheorem{proposition}{Proposition}[section]
\newtheorem{corollary}{Corollary}[section]
\newtheorem{theorem}{Theorem}[section]
\newtheorem{remark}{Remark}[section]

\newcommand{\s}{\vspace{2ex}}
\newcommand{\n}{\noindent}
\newcommand{\e}{\varepsilon}
\newcommand{\R}{\mathbb{R}}
\newcommand{\Z}{\mathbb{Z}}

\newcommand{\E}{\mathbb{E}}
\newcommand{\PP}{\mathbb{P}}

\title{Lyapunov exponents and shear-induced chaos for a Hopf bifurcation with additive noise}

\author{Peter H. Baxendale} 

\AtEndDocument{\bigskip{\footnotesize
  \noindent
  \textsc{Peter H. Baxendale. Department of Mathematics, University of Southern California, CA 90089, USA.}
  \textit{E-mail address}: \texttt{\href{mailto:baxendal@usc.edu}{baxendal@usc.edu}}
 }}

\begin{document}

\maketitle

\begin{abstract}

This paper considers the effect of additive white noise on the normal form for the supercritical Hopf bifurcation in 2 dimensions. The main results involve the asymptotic behavior of the top Lyapunov exponent $\lambda$ associated with this random dynamical system as one or more of the parameters in the system tend to 0 or $\infty$. This enables the construction of a bifurcation diagram in parameter space showing stable regions where $\lambda <0$ (implying synchronization) and unstable regions where $\lambda > 0$ (implying chaotic behavior). The value of $\lambda$ depends strongly on the shearing effect of the twist factor $b/a$ of the deterministic Hopf bifurcation.  If $b/a$ is sufficiently small then $\lambda <0$ regardless of all the other parameters in the system. But when all the parameters except $b$ are fixed then $\lambda$ grows like a positive multiple of $b^{2/3}$ as $b \to \infty$.

\end{abstract}

\s

\n Keywords: random dynamical system, Lyapunov exponent, stochastic Hopf bifurcation, shear-induced chaos.

\section{Introduction} \label{sec intro}

Consider the stochastic differential equation (SDE) in $\R^2$:
 \begin{equation} \label{X}
  dX_t  = \left( \begin{bmatrix} \mu  & -\omega\\
  \omega & \mu \end{bmatrix}X_t + \|X_t\|^2\begin{bmatrix} -a & -b \\ b & -a \end{bmatrix}X_t\right)dt + \sigma \begin{bmatrix} dW_t^1 \\ dW_t^2 \end{bmatrix} \end{equation}
for constants $\mu,\omega,a,b$ and $\sigma$, where $\{W_t^i: t \ge 0\}$ for $i = 1,2$ are independent standard Brownian motions.  Assume $\sigma > 0$, and to ensure the system is recurrent assume $a >0$.  The other constants $\mu, \omega, b$ can be arbitrary real numbers.  Since the noise is additive the It\^o and Stratonovich versions of \eqref{X} agree.

When $\sigma = 0$ the equation reduces to the normal form for the supercritical Hopf bifurcation with fixed point $0 \in \R^2$ and bifurcation parameter $\mu$.  See for example Guckenheimer and Holmes 
\cite{GH90}.  If $\mu < 0$ then the fixed point $0 \in \R^2$ is globally stable, and if $\mu >0$ then $0$ is unstable and there is stable limit cycle centered at $0$ with radius $\sqrt{\mu/a}$.   The parameter $a$ determines the local stability of the limit cycle and the parameter $b$ determines the shearing effect of the system near the limit cycle.      

%We will study the effect of additive noise on this system, and henceforth assume $\sigma > 0$. 
The stochastic system with $\sigma >0$ has been studied from a theoretical viewpoint in DeVille, Namachchivaya and Rapti \cite{DNR},  Doan, Engel, Lamb and Rasmussen \cite{DELR},  Chemnitz and Engel \cite{CE23} and Breden and Engel \cite{BE23}.  See also  Wieczorek \cite{Wie} in relation to stochastic bifurcation of noise-driven lasers, and Dijkstra, Frankcombe and von der Heydt \cite{DFH08} and Tantet, Chekroun, Dijkstra and Neelin \cite{CTDN-II} in relation to models in oceanography.

We consider how the stability of the system depends on the parameters $\mu$, $\omega$, $a$, $b$ and $\sigma$.  That is, we wish to describe the stochastic bifurcation scenario.  In one dimension,  Crauel and Flandoli \cite{CF98} show that ``Additive noise destroys a pitchfork bifurcation.''   
The two-dimensional nature of the Hopf bifurcation makes the stability analysis more complicated.  Monotonicity arguments are not available.  
Instead, the extra dimension allows the possibility of shearing vector fields, and it turns out that the size of the parameter $b$ controlling the strength of the shearing is important in determining the stability of the system.

The term ``shear induced chaos'' refers to the tendency of shear vector fields to cause chaotic behavior, and this has been well-studied in both the deterministic and stochastic settings.  
Of particular relevance is the study of periodically and stochastically ``kicked'' limit cycles, see Wang and Young \cite{WY03} and Lin and Young \cite{LY08}, and Lu and Wang and Young \cite{LWY13}. 
 When $\mu > 0$ the deterministic system underlying \eqref{X} has a stable limit cycle, and the white noise acts as random ``kicking'' in the two coordinate directions.  
 Engel, Lamb and Rasmussen \cite{ELR} consider a special case of stochastic kicking of a limit cycle, and their formula for the associated Lyapunov exponent plays an important role in the recent result of Chemnitz and Engel \cite{CE23} for the system \eqref{X}.

 With additive white noise the one-point motion $\{X_t: t \ge 0\}$ is positive recurrent on $\R^2$, with invariant density with respect to Lebesgue measure of the form 
       $$
       \frac{1}{Z} \exp\left(- \frac{a}{2\sigma^2}\bigl(\|X\|^2-\mu/a\bigr)^2 \right)
       $$ 
for some normalizing constant $Z$.  The shape of the invariant density varies with $\mu$, and there is a P-bifurcation (see Arnold \cite{Arn98}) at $\mu = 0$, but there is no essential qualitative change in the dynamical behavior of the one-point motion.  In this setting it is more interesting to consider the associated random dynamical system, and in particular the issue of stability along trajectories.  
  Suppose that $\{X_t: t\ge 0\}$ and $\{Y_t: t \ge 0\}$ are both strong solutions of \eqref{X}
with distinct initial conditions $X_0$ and $Y_0$.  It is a non-trivial question to determine whether or not
$\|Y_t-X_t\| \to 0$ almost as $t\to \infty$.  If $\|Y_t-X_t\| \to 0$ almost surely for all distinct $X_0$ and $Y_0$ in $\R^2$ we say that the system \eqref{X} has stability along trajectories.  In applications involving oscillators (see for example \cite{Wie,LSBY09,PRK}) this is often referred to as synchronization, or phase locking. 

In this paper we study the closely related concept of linearized stability along trajectories.  Linearizing \eqref{X} along a trajectory $\{X_t: t \ge 0\}$ gives the process $\{U_t: t \ge 0\}$ satisfying 
\begin{equation} \label{U}
   dU_t = \left( \begin{bmatrix} \mu  & -\omega\\
  \omega & \mu \end{bmatrix}U_t + \|X_t\|^2\begin{bmatrix} -a & -b \\ b & -a \end{bmatrix}U_t +  2 \langle X_t,U_t \rangle \begin{bmatrix} -a & -b \\ b & -a \end{bmatrix}X_t \right)dt.
  \end{equation} 
Here $U_t$ should be regarded as a vector at the point $X_t$.  Intuitively, if $U_0 = Y_0-X_0$ and $\|Y_0-X_0\|$ is small then the processes $\{U_t: t \ge 0\}$ and $\{Y_t-X_t: t \ge 0\}$ should evolve in a
similar manner as long as $\|Y_t-X_t\|$ remains small.  The long time rate of exponential growth or decay of the process $\{U_t: t \ge 0\}$ is given by the (top) Lyapunov exponent   
   \begin{equation} \label{lam def}
    \lambda = \lambda(\mu,\omega,a,b,\sigma) = \lim_{t \to \infty} \frac{1}{t} \log \|U_t\|.
    \end{equation}
For the system (\ref{X},\ref{U}) standard arguments show that the limit in \eqref{lam def} exists almost surely and takes the same value for all initial points $X_0 \in \R^2$ and all non-zero initial vectors $U_0 \in \R^2$.

Doan, Engel, Lamb and Rasmussen \cite{DELR} show that the system \eqref{X} gives rise to a random dynamical system in the sense of Arnold \cite{Arn98}, and they discuss the relationship between the sign of $\lambda$ and the nature of the random attractor associated with the random dynamical system.  In particular they show, applying a result of Flandoli, Gess and Scheutzow \cite{FGS}, that $\lambda <0$ for \eqref{X} implies synchronization, so that the random attractor is almost surely a singleton.  In this paper we concentrate on the evaluation of $\lambda$ and especially the determination of regions in parameter space $(\mu,\omega,a,b,\sigma)$ where $\lambda$ is negative (indicating stability along trajectories) or positive (indicating sensitive dependence on initial conditions, and chaos).

 \begin{remark} \label{rem mult bif} The bifurcation scenario is very different if the additive noise in \eqref{X} is replaced by multiplicative noise, so that the deterministic fixed point $0$ remains fixed under the perturbation.  For example, consider the effect of varying the parameter $\mu$ in the Stratonovich SDE with linear white noise:
   \begin{equation} \label{X mult}
   dX_t = \begin{bmatrix} \mu & -\omega \\ \omega & \mu \end{bmatrix}X_tdt 
      + \|X_t\|^2 \begin{bmatrix} -a & -b \\ b & -a \end{bmatrix}X_tdt + \sum_{\alpha = 1}^n A_\alpha X_t \circ  dW_t^\alpha
   \end{equation}
for some given matrices $A_1, \ldots, A_n$ (and assume at least one of the $A_\alpha$ is not skew-symmetric).  Linearizing \eqref{X mult} at 0 gives
    \begin{equation} \label{U mult}
   dU_t = \begin{bmatrix} \mu & -\omega \\ \omega & \mu \end{bmatrix}U_tdt 
      + \sum_{\alpha = 1}^n  A_\alpha U_t \circ dW_t^\alpha.
   \end{equation}
Let $-\mu_c$ denote the top Lyapunov exponent for the SDE
    $$
   dV_t = \begin{bmatrix} 0 & -\omega \\ \omega & 0 \end{bmatrix}V_tdt 
      + \sum_{\alpha = 1}^n  A_\alpha V_t \circ dW_t^\alpha
  $$ 
satisfied by $V_t = e^{-\mu t}U_t$.  Then the top Lyapunov exponent for \eqref{U mult} is $\lambda = \mu - \mu_c$.  If $\mu < \mu_c$ then $0$ is a stable fixed point and $X_t \to 0$ almost surely for all $X_0 \in \R^2$.  But if $\mu > \mu_c$ then $0$ is an unstable fixed point and $X_t \to \nu$ in distribution for all $X_0 \neq 0$, where $\nu$ is an invariant probability measure for $\{X(t): t \ge 0\}$ on $\R^2 \setminus \{0\}$.  Moreover $\nu$ is close to 0 in the sense that for each $p >0$ there exists $C_p$ such that $\int_{\R^2 \setminus \{0\}} \|X\|^p d\nu(X) \sim C_p(\mu-\mu_c)$ as $\mu \searrow \mu_c$.  This is an application of Baxendale \cite[Theorem 2.13]{Bax94}, which applies also to more general multiplicative noise systems in arbitrary dimensions.  The issue of stability along trajectories for the system in $\R^2 \setminus \{0\}$ is addressed in Baxendale \cite{Bax97}.

   \end{remark}

\subsection{Evaluation of $\lambda$}  \label{sec eval}
The rotational invariance of the system \eqref{X} simplifies the computation of $\lambda$ considerably.  The Furstenberg-Khasminskii formula (see Arnold \cite{Arn98}) gives
   \begin{equation} \label{lam khas}
   \lambda = \lim_{t \to \infty} \frac{1}{t} \int_0^t Q(r_s, \psi_s)ds = \int_M Q(r,\psi)d\widetilde{\nu}(r,\psi)
   \end{equation}    
where
    \begin{align} 
    Q(r,\psi)  & =  \mu-2ar^2 +r^2\bigl( b \sin 2\psi- a \cos 2\psi \bigr) \nonumber \\
    & = \mu-ar^2 +2r^2 \cos \psi \bigl( b \sin \psi- a \cos  \psi \bigr) \label{Q}
    \end{align}
and $\widetilde{\nu}$ is the unique invariant probability measure on $M \equiv (0,\infty) \times \R/(\pi \Z)$ for the diffusion $\{(r_t,\psi_t): t \ge 0\}$ with generator ${\cal L}$, say, given by
    \begin{align} \label{r}
   dr_t & = \bigl(\mu r_t - ar_t^3+\frac{\sigma^2}{2r_t}\bigr)dt + \sigma dW_t^r,\\
    d\psi_t & =  2r_t^2 \cos \psi_t \bigl(b \cos \psi_t+a \sin \psi_t\bigr)dt-\frac{\sigma}{r_t}dW_t^\phi, \label{psi}
   \end{align} 
where $\{W_t^r: t \ge 0\}$ and $\{W_t^\phi: t \ge 0\}$ are independent standard Brownian motions.  The derivation of these formulas is given in Section \ref{sec FK}.  A version with a slightly different parametrization for $\psi$ is given in \cite[Section 3.1]{DNR}.

 The rotation parameter $\omega$ does not appear in the equations (\ref{Q},\ref{r},\ref{psi}), and so $\lambda$ does not depend on $\omega$.  The underlying reason for the lack of dependence on $\omega$ is that the system \eqref{X} is equivariant in distribution under the group $SO(2)$ of rigid rotations of $\R^2$.  Henceforth we write 
    $$
    \lambda = \lambda(\mu,a,b,\sigma) \quad \mbox{where} \quad (\mu,a,b,\sigma) \in \R \times (0,\infty) \times \R \times (0,\infty).
    $$
Also the function $Q(r,\psi)$ is invariant under the transformation $(b,\psi) \mapsto (-b,-\psi)$ and the system (\ref{r},\ref{psi}) is invariant in distribution under the same transformation.  It follows that 
  \begin{equation} \label{b}
  \lambda(\mu,a,b,\sigma) = \lambda(\mu,a,-b,\sigma).
  \end{equation}
In this case the underlying reason is that the system \eqref{X} is equivariant in distribution under the reflection $\begin{bmatrix} x_1\\x_2\end{bmatrix} \to \begin{bmatrix} x_1\\-x_2\end{bmatrix}$.  Because of this symmetry we will assume without loss of generality that $b \ge 0$.  We note that \cite{DELR} and \cite{CE23} use $-b$ where \cite{DNR} and this paper use $b$; happily, because of \eqref{b} this is not a source of confusion.  %e and consider $\lambda(\mu,a,b,\sigma)$ for $(\mu,a,b,\sigma) \in \R \times (0,\infty) \times [0,\infty) \times (0,\infty)$.  
  
\s

There is no explicit formula for the density of $\widetilde{\nu}$, so it is hard to evaluate $\lambda$ exactly.   However there is an explicit formula \eqref{rho} for the density for the $r$ marginal $\nu$, say, of $\widetilde{\nu}$ and the methods in this paper will make much use of this fact.  
% {\bf This next is probably a repeat of earlier remarks.}  
The results in this paper will describe the asymptotic behavior of $\lambda$ as one or more of the parameters tends to 0 or infinity.  These will lead to results giving regions in parameter space for which $\lambda > 0$ or $\lambda < 0$.

\subsection{Structure of the paper}

  After a brief review of earlier results in Section \ref{sec prev}, we present the main results in Section \ref{sec results}.  Theorem \ref{thm binf} shows that $\lambda(\mu,a,b,\sigma)$ grows asymptotically like $b^{2/3}$ as $b \to \infty$, so that $\lambda > 0$ for sufficiently large $b$.   A uniform version of this result yields a confirmation of Conjecture D of Doan et al \cite{DELR} concerning the positivity of $\lambda$, see Corollary \ref{cor binf unif2}.  
  On the other hand the small noise asymptotic formula in Theorem \ref{thm sigma0}, together with the scaling result Proposition \ref{prop scale}, leads to the existence of an absolute constant $k_0 > 0$ such that $\lambda(\mu,a,b,\sigma) <0$ whenever $0 \le b/a \le k_0$, see Corollary \ref{cor sigma0 2}.
   Sections \ref{sec stocbif} and \ref{sec stab diag} contain the consequences of these results for the bifurcation scenario and the corresponding stability diagram.   

Section \ref{sec FK} gives the derivation of the Furstenberg-Khasminskii formula \eqref{lam khas} and contains a brief discussion of the adjoint method, introduced by Arnold, Papanicolaou and Wihstutz \cite{APW86}, which is adapted here in several different ways to prove the main results.  In Section \ref{sec h} we give some heuristic arguments in favor of these results.  This a non-rigorous discussion, but in the case of Theorems \ref{thm binf} and \ref{thm CE unif} especially it motivates the method of rigorous proof later.   The proofs of Theorems \ref{thm binf} to \ref{thm CE unif} are given in Sections \ref{sec binf proof} to \ref{sec sigma0 proof}.  Finally the Appendix contains some auxiliary results about the invariant probability measure $\nu$ for the process $\{r_t: t \ge 0\}$ as well as a technical result Proposition \ref{prop-bg} allowing the implementation of the adjoint method on a non-compact state space.

\section{Previous results} \label{sec prev}

\subsection{An upper bound on $\lambda$} \label{sec prev upper}

This section contains results in DeVille, Namachchivaya and Rapti \cite{DNR} with extensions by Doan, Engel, Lamb and Rasmussen \cite{DELR}.  The upper bound     
    \begin{equation} \label{Qbound}  Q(r,\psi) \le \mu +(\sqrt{a^2+b^2}-2a)r^2
   \end{equation}
and the explicit formula for $\int_{\R_+} r^2d\nu(r)$ given in \cite[eqn (26)]{DNR}, see \eqref{int r2} in the Appendix, together imply
   \begin{align}  \nonumber 
   \lambda(\mu,a,b,\sigma)  =\int_M Q(r,\psi)d\widetilde{\nu}(r,\psi) & < \mu + (\sqrt{a^2+b^2}-2a) \int_{\R_+} r^2 d\nu(r)\\
   &  = \sqrt{2a \sigma^2} J\left(\frac{\mu}{\sqrt{2a\sigma^2}},\frac{b}{a}\right) \label{lam J}
   \end{align}
where 
   \begin{equation} \label{J}
   J(z,b) = z+ (-2+\sqrt{1+b^2})\left(z+ \frac{1}{\sqrt{\pi}} \frac{ e^{-z^2}}{  {\rm erfc}(-z)}\right).
   \end{equation}
The condition $J(z,b) \le 0$ is equivalent to the condition $|b| \le \widehat{J}(z)$ where 
     \begin{equation} \label{Jhat}
     \widehat{J}(z) =  \sqrt{\dfrac{1}{1+ \sqrt{\pi} ze^{z^2}{\rm erfc}(-z)}\left(\dfrac{1}{1+ \sqrt{\pi} ze^{z^2}{\rm erfc}(-z)}+2 \right)} > 0.
     \end{equation}
It follows from elementary analysis that $\widehat{J}(z)$ is strictly decreasing with $\widehat{J}(z) \sim 2 z^2$ as $z \to -\infty$ and $\widehat{J}(z) \sim 1/\sqrt{\sqrt{\pi}z e^{z^2}}$ as $z \to \infty$.  The graph of $\mu \mapsto \widehat{J}(\mu/\sqrt{2})$ appears as part of Figure \ref{fig stab diag}.  The inequality \eqref{lam J} gives the following result.

\begin{proposition} \label{prop Jhat} If $0 \le b \le a \widehat{J}\left(\dfrac{\mu}{\sqrt{2a \sigma^2}}\right)$ then $\lambda(\mu,a,b,\sigma) < 0$.  In particular, given $\mu \in \R$, $a >0$ and $\sigma > 0$ there exists $b_0 > 0$ such that $\lambda(\mu,a,b,\sigma) < 0$ whenever $0 \le b \le b_0$.%if $b \le \widehat{J}(\mu/\sqrt{2})$ then $\lambda(\mu,1,b,1) < 0$.  In particular
\end{proposition} 
 
This includes the result $\lambda(\mu,a,0,\sigma) < 0$ given by Flandoli, Gess and Scheutzow \cite{FGS}.

\begin{remark} \label{rem norm}  The equation (10) in \cite{DNR} for the normalization constant $Z$ for the stationary density of $\{(x_t,y_t): t \ge 0\}$ in $\R^2$ is incorrect by a factor $2 \pi e^{\mu^2/2a\sigma^2}$. 
This incorrect value for the normalization constant is not used in \cite{DNR}, but it is quoted and used in  \cite{DELR}.  In the notation of \cite{DELR}, using their $\alpha$ in place of our $\mu$, the normalization constant $K_{a,\alpha,\sigma}$ in \cite[eqn (2.3)]{DELR} should be  
    $$
    K_{a,\alpha,\sigma} = \frac{2 \sqrt{2a}}{\sigma \sqrt{\pi} \,{\rm erfc}(-\alpha/\sqrt{2a\sigma^2})}
   \cdot\left(2 \pi e^{\alpha^2/2a\sigma^2}\right)^{-1}  = \frac{\sqrt{2a}}{\sigma \pi^{3/2}e^{-\alpha^2/2a\sigma^2} \,{\rm erfc}(-\alpha/\sqrt{2a\sigma^2})}
    $$
The statement and proof of Theorem C in \cite{DELR} are correct, except that the wrong value for $K_{a,\alpha,\sigma}$ is used.  As soon as the correct value for $K_{a,\alpha,\sigma}$ is used, then the condition in Theorem C of \cite{DELR} becomes exactly the same as in Proposition \ref{prop Jhat} above.  As  $\alpha \to 0$ and as $\sigma \to \infty$ the ratio of the incorrect $K$ and the correct $K$ converges to $2 \pi$, so that the statements in Remark 2.1 of \cite{DELR} remain valid. 
 
 \end{remark}
 
\begin{remark} \label{rem unif}  The upper bound \eqref{Qbound} and the corresponding lower bound $Q(r,\psi) \ge \mu-(\sqrt{a^2+b^2}+2a)r^2$ are both sharp for $\psi \in \R/\pi\Z$.  It follows that $Q(r,\psi)$ is essentially unbounded below for the invariant measure $\widetilde{\nu}$, and is essentially bounded above if and only if $b \le \sqrt{3}a$.  If $b \le \sqrt{3}a$ then the upper bound $Q(r,\psi) \le \mu$ is sharp.  These facts about $Q(r,\psi)$, together with the expression \eqref{norm} for $\log\|U_t\|$ in terms of $Q$, are very closely related to, and help illuminate, the discussion of dichotomy spectrum, bounds on finite time Lyapunov exponents and issues of uniform or non-uniform attractivity in Theorems E and F of \cite{DELR}.  

\end{remark}  
 
\subsection{Result of Chemnitz and Engel} \label{sec prev CE}

Let $\Psi(\zeta)$ denote the top Lyapunov exponent for the SDE
     \begin{equation} \label{YB}
            dY_t = \begin{bmatrix} -1  & 0 \\ \zeta^{1/3} & 0 \end{bmatrix}Y_t dt + \begin{bmatrix} 0 & \zeta^{1/6} \\ 0 & 0 \end{bmatrix} Y_t dW_t. 
      \end{equation}      
Imkeller and Lederer \cite{IL01} show that
    \begin{equation} \label{Psi}
   \Psi(\zeta) = \frac{1}{2}\left( \dfrac{ \int_0^\infty u^{1/2} \exp\left(-\frac{1}{\zeta}\left(\frac{1}{6}u^3-\frac{1}{2}u \right)\right)du}{ \int_0^\infty u^{-1/2} \exp\left(-\frac{1}{\zeta}\left(\frac{1}{6}u^3-\frac{1}{2}u \right)\right)du}-1\right).
   \end{equation} 
Engel, Lamb and Rasmussen \cite{ELR} show there exists $c_\ast > 0$ such that $\Psi(\zeta) < 0$ for $0 < \zeta < c_\ast$ and $\Psi(\zeta) > 0$ for $\zeta > c_\ast$.  Numerically, $c_\ast \approx (0.2823)^{-1} = 3.543$.   

\begin{theorem} \label{thm CE} {\bf (Chemnitz and Engel \cite{CE23})}  Fix $\mu > 0$, $a >0$, $b \ge 0$ and $\sigma >0$.  Then
  $$
  \lim_{\e \to 0} \lambda(\mu,a,b/\e,\e\sigma) = 2 \mu \Psi\left(\frac{b^2 \sigma^2}{2 \mu^2 a}\right).
  $$
In particular,  if $b^2 \sigma^2 > 2c_\ast \mu^2 a$ then $\lambda(\mu,a,b/\e,\e \sigma) > 0$ for all sufficiently small $\e$.
\end{theorem}
   
This gives the first proof that $\lambda > 0$ is possible in this system.

\section{Main results} \label{sec results}

\subsection{Scaling}   \label{sec scaling}  We have already seen that $\lambda$ does not depend on $\omega$ or the sign of $b$.  The next result reduces the effective parameter space to 2 dimensions.

\begin{proposition} \label{prop scale}  For all $A >0$ and $B >0$ we have
 \begin{equation} \label{scale} \lambda(\mu,a,b,\sigma) = \frac{1}{A}\lambda\left(A\mu,\frac{Aa}{B^2},\frac{Ab}{B^2},\sqrt{A}B\sigma\right). 
      \end{equation} 
In particular 
      \begin{equation} \label{scale2} \lambda(\mu,a,b,\sigma) = \sigma\sqrt{a} \,\lambda\left(\frac{\mu}{\sigma \sqrt{a}},1,\frac{b}{a},1\right)
   \end{equation} 
and if $\mu \neq 0$ then  
    \begin{equation} \label{scale3} \lambda(\mu,a,b,\sigma) = |\mu| \lambda\left(\frac{\mu}{|\mu|},a,b,\frac{\sigma }{|\mu|}\right). 
   \end{equation} 
\end{proposition} 

\n{\bf Proof.}    
Rescaling the system \eqref{X} in time by a factor $A$ gives an equivalent system with parameters $(A\mu,Aa,Ab,\sqrt{A}\sigma)$, so that $A\lambda(\mu,a,b,\sigma) = \lambda(A\mu,Aa,Ab,\sqrt{A}\sigma)$.   Also, rescaling the system \eqref{X} in space by a factor $B$ gives an equivalent system with parameters $(\mu,a/B^2,b/B^2,B\sigma)$, so that $\lambda(\mu,a,b,\sigma) = \lambda(\mu,a/B^2,b/B^2,B\sigma)$.  Together we have \eqref{scale}.   Taking $A = \sigma \sqrt{a}$ and $B = \sqrt{Aa}$ gives \eqref{scale2}, and taking $A = 1/|\mu|$ and $B = \sqrt{A}$ gives \eqref{scale3}.   \qed

\s

It follows from \eqref{scale2} that the sign of $\lambda$ depends only on the dimensionless ratios $b/a$ and $\mu/(\sigma \sqrt{a})$.  The ratio $b/a$ is the ``twist factor'' or ``twist number'' of the deterministic Hopf bifurcation obtained by putting $\sigma = 0$ in \eqref{X}, see for example \cite{WY03,LWY13,CTDN-II}. When $\mu >0$ it determines the angle of transversality as the deterministic isochrons intersect the limit cycle of radius $\sqrt{\mu/a}$.  The interpretation of $\mu/(\sigma \sqrt{a})$ is less clear.  Adopting the argument in \cite{CE23} for the case of positive $\mu > 0$ and small noise intensity $\sigma$, we have {\rm radius} $ = \sqrt{\mu/a}$, {\rm contraction} $= 2\mu$ and {\rm effective noise = noise/radius} $= \sigma \sqrt{a/\mu}$.  Then 
    $$
    \frac{\rm contraction}{({\rm effective\,\, noise})^2} = \frac{2\mu}{\sigma^2 a/\mu} = 2\left(\frac{\mu}{\sigma\sqrt{a}}\right)^2.
    $$
But this interpretation is not valid when $\mu$ is negative or when $\sigma$ is not small.

%For the stability diagram it is enough to consider $\lambda(\mu,1,b,1)$ for $(\mu,b)$ in the half space $\R \times [0,\infty)$.  The equation \eqref{scale3} will be useful for converting small $\sigma$ results into large $\mu$ results. 

%\begin{remark}  The effect of the scaling result \eqref{scale2} can be seen in the numerical results of Wieczorek \cite{Wie}.  The straight lines in Figure 8 of \cite{Wie} and the threefold repetition of the same curve in Figure 9 of \cite{Wie} are direct consequences of the scaling.
%\end{remark}

\subsection{Results showing $\lambda > 0$} \label{sec pos}

Let $\gamma_0$ denote the top Lyapunov exponent for the 2-dimensional linear SDE
       \begin{equation} \label{YA}
 dY_t = \begin{bmatrix} 0 & 0 \\1 & 0 \end{bmatrix}Y_tdt
   +\begin{bmatrix} 0 & 1 \\0 & 0 \end{bmatrix}Y_tdW_t.
   \end{equation}
Ariaratnam and Xie \cite{AX90} show that $\gamma_0 = \dfrac{\pi}{2^{1/3} 3^{1/6}[\Gamma(1/3)]^2}\approx 0.29$.
 Write $\nu = \nu_{\mu,a,\sigma}$ to show the dependence of the $\{r_t: t \ge 0\}$ invariant probability measure $\nu$ on the parameters $\mu \in \R$, $a > 0$ and $\sigma >0$.   
    
\begin{theorem} \label{thm binf}   Given $\mu \in \R$, $a > 0$ and $\sigma >0$,     
   \begin{equation} \label{binf asymp}
   \lambda(\mu,a,b,\sigma) \sim (2 b \sigma)^{2/3} \gamma_0 \int_{\R_+} r^{2/3} d\nu_{\mu,a,\sigma}(r) \quad \mbox{ as } b \to \infty.
    \end{equation}
\end{theorem}

The proof of Theorem \ref{thm binf} is given in Section \ref{sec binf proof}.  Note the $b^{2/3}$ growth rate.  A similar fractional power shows up for small random perturbations of 2 dimensional Hamiltonian systems, see Baxendale and Goukasian \cite{BG02} and the conjecture in Arnold, Imkeller and Namachchivaya \cite{AIN04}.  %The reason is the same in all these cases; the system can be modelled as 
%The setting here is closely related to that in these earlier papers; 
Just as in these earlier papers, the result crucially involves the change of basis for linear nilpotent It\^{o} systems introduced in the seminal paper of Pinsky and Wihstutz \cite{PW88}.  More details are given in the heuristic argument in Section \ref{sec h binf}.  

Theorem \ref{thm binf} implies that given $\mu \in \R$, $a > 0$ and $\sigma >0$ there exists $b_0 >0$ such that $\lambda(\mu,a,b,\sigma) >0$ whenever $b \ge b_0$.  This implies shear-induced chaos for sufficiently large shear parameter $b$.  Much of the discussion of shear induced chaos is in the setting of perturbations of limit cycles, which corresponds here to the case of $\mu > 0$.  Note that Theorem \ref{thm binf} includes the case $\mu < 0$ for which there is no limit cycle to be perturbed by noise.

The next result follows from a detailed analysis of the convergence rate in \eqref{binf asymp}.

%Theorem \ref{thm binf} implies for fixed $\mu$, $a$ and $\sigma$ then $\lambda(\mu,a,b,\sigma) > 0$ for all sufficiently large $b$.  The next result is based on uniform estimates on the rate of convergence in Theorem \ref{thm binf} and the positivity of the right side in \eqref{binf asymp}.  

\begin{theorem} \label{thm binf unif}  Suppose $P \subset \R \times (0,\infty) \times (0,\infty)$ is bounded and that $\int_{\R_+} r^2 d\nu_{(\mu,a,\sigma)}(r)$ and $\int_{\R_+} r^{-2/3} d\nu_{(\mu,a,\sigma)}(r)$ are bounded above for $(\mu,a,\sigma) \in P$. 

(i) For all $\e >0$ there exists $k < \infty$ such that
    \begin{equation} \label{binf unif}
       \left| \frac{\lambda(\mu,a,b,\sigma)}{(2b\sigma)^{2/3}} - \gamma_0 \int_{\R_+} r^{2/3}d\nu_{(\mu,a,\sigma)}(r)\right| <\e
       \end{equation}
whenever $(\mu,a,\sigma) \in P$ and $b\sigma \ge k$. 

(ii) If also $\int_{\R_+} r^{2/3} d\nu_{(\mu,a,\sigma)}(r)$
is bounded below away from 0 for $(\mu,a,\sigma) \in P$ then there exists $k  <\infty$ such that $\lambda(\mu,a,b,\sigma) > 0$ whenever $(\mu,a,\sigma) \in P$ and $b\sigma \ge k$.  \end{theorem}
   
Three different choices for the set $P$ give the following three consequences.    

\begin{corollary} \label{cor binf unif}  (i)  Fix $a > 0$ and $\sigma > 0$.  Given $\mu_0 > 0$ there exists $k_1 < \infty$ such that $\lambda(\mu,a,b,\sigma) >0$ if $-\mu_0 \le \mu \le \mu_0$ and $b \ge k_1$.

(ii) Fix $\mu > 0$ and $a > 0$.  Given $\sigma_0 > 0$ there exists $k_2 < \infty$ such that $\lambda(\mu,a,b,\sigma) >0$ if $0 < \sigma \le \sigma_0$ and $b\sigma  \ge k_2$.

(iii) %Fix $\mu < 0$ and $\sigma > 0$.  Given $a_0 > 0$ there exists $k_3$ such that $\lambda(\mu,a,b,\sigma) >0$ if $0 < a \le a_0$ and $b \ge k_3$.
   Fix $\mu <0$ and $a > 0$.  Given $\sigma_0 > 0$ there exists $k_3< \infty$ such that $\lambda(\mu,a,b,\sigma) > 0$ if $0 <\sigma \le \sigma_0$ and $b\sigma^2 \ge k_3$.
   \end{corollary}

The proofs of Theorem \ref{thm binf unif} and Corollary \ref{cor binf unif} are given in Section \ref{sec proof unif}.  %the verification that each In each case the verification that $P$ satisfies the conditions of Theorem \ref{thm binf unif} is given LATER.%\begin{remark}
  Notice that the lower bounds for $b$ involve different powers of $\sigma$ in each part of Corollary \ref{cor binf unif}.  
  In (i) the noise parameter $\sigma$ is fixed and so a lower bound on $b\sigma$ is equivalent to a lower bound on $b$.  
  In (iii) since $\mu < 0$ then $\nu_{(\mu,a,\sigma)}$ converges to the point mass at 0 as $\sigma \to 0$, and so the natural parameter set $P = \{(\mu,a,\sigma): 0 < \sigma \le \sigma_0\}$ does not satisfy the assumptions for either part of Theorem \ref{thm binf unif}.  
 Taking $B = 1/\sigma$ in the proof of Proposition \ref{prop scale} we see that $\{r_t/\sigma: t \ge 0\}$ has invariant probability measure $\nu_{(\mu,a\sigma^2,1)}$ so that $\sigma^{-2/3}\int_{\R_+} r^{2/3}d\nu_{(\mu,a,\sigma)}(r) =  \int_{\R_+}r^{2/3}d\nu_{(\mu,a\sigma^2, 1)}(r) $ which is bounded away from 0 for $0 <\sigma \le \sigma_0$.   This suggests that a different power of $\sigma$ will be needed.  
  The proof of (iii) uses the parameter set $P= \{(\mu,a\sigma^2,1): 0 < \sigma \le \sigma_0\}$ together with the scaling result $\lambda(\mu,a,b,\sigma) = \lambda(\mu,a\sigma^2, b\sigma^2,1)$ obtained by taking $B = 1/\sigma$ in \eqref{scale}.  This explains why the lower bound is on $b\sigma^2$ rather than on $b\sigma$.
%    \end{remark}

\s

Applying the scaling results in Proposition \ref{prop scale} we get the following confirmation of Conjecture D of Doan, Engel, Lamb and Rasmussen \cite{DELR}.

\begin{corollary} \label{cor binf unif2}  For each $c >0$ there exist constants $k_1$, $k_2$ and $k_3$ such that
$\lambda(\mu,a,b,\sigma) > 0$ whenever $b \ge aC\left(\dfrac{\mu}{\sqrt{a}\sigma}\right
)$ with the function
     $$
     C(x) = \begin{cases} k_3 x^2 & \mbox{ if }x < -c, \\
                             k_1 & \mbox{ if } -c \le x \le c, \\
                              k_2 x & \mbox{ if } x > c.
                              \end{cases}
                            $$
\end{corollary}

\n{\bf Proof.}  By \eqref{scale2} it suffices to prove the result for the sign of $\lambda(\mu,1,b,1)$.  Fix $c >0$.  Apply Corollary \ref{cor binf unif}(i) with $a= \sigma = 1$ and $\mu_0 = c$ to obtain $k_1$.  Then $\lambda(\mu,1,b,1) >0$ if $-c\le \mu \le c$ and $b \ge k_1$.  Next apply  Corollary \ref{cor binf unif}(ii) with $\mu = a= 1$ and $\sigma_0 =1/c$ to obtain $k_2$ so that $\lambda(1,1,b,\sigma) > 0 $ if $0 <\sigma \le 1/c$ and $b\sigma \ge k_2$.  Using \eqref{scale3} we have $\lambda(\mu,1,b,1) = \lambda(1,1,b,1/\mu) >0$ if $\mu \ge c$ and $b \ge k_2 \mu$.   Finally apply Corollary \ref{cor binf unif}(iii) with $\mu = -1$, $a= 1$ and $\sigma_0 = 1/c$ to obtain $k_3$ so that $\lambda(-1,1,b,\sigma) > 0 $ if $0 <\sigma \le 1/c$ and $b\sigma^2 \ge k_3$.  Using \eqref{scale3} with $\mu < 0$ we have $\lambda(\mu,1,b,1) = \lambda(-1,1,b,1/|\mu|) >0$ if $|\mu| \ge c$, that is $\mu \le -c$, and $b \ge k_3 \mu^2$.  \qed

\subsection{Results showing $\lambda <0$} \label{sec neg}

\begin{theorem} \label{thm sigma0} Fix $\mu > 0$ and $a > 0$.
For all $\sigma_0 > 0$ and $k >0$ there exists $K < \infty$ such that 
    $$
    \left|\lambda(\mu,a,b,\sigma) + \frac{(a^2+b^2)\sigma^2}{2 \mu a}\right| \le K \sigma^4
    $$
whenever $0 \le b  \le k$ and $0 < \sigma \le \sigma_0$.
\end{theorem}

The proof of Theorem \ref{thm sigma0} is given in Section \ref{sec sigma0 proof}.

\begin{remark} DeVille, Namachchivaya and Rapti \cite[Prop 4.1(2)]{DNR} claim $\lambda(\mu,a,b,\sigma) = -3a\sigma^2/(4 \mu) + {\cal O}(\sigma^4)$ as $\sigma \to 0$ when $\mu > 0$.  Their calculations involve an approximation for $\widetilde{\nu}$ as the product of independent Gaussian measures.  However,  this approximation is incapable of correctly capturing all the order $\sigma^2$ terms in the evaluation of the Furstenberg-Khasminskii formula \eqref{lam khas}.  See for example the difference at order $\sigma^2$ in the evaluation of $\E[r^2]$ using the exact density $\nu_r(r)$ given in \cite[eqn (21)]{DNR}, see \cite[eqn (26)]{DNR}, and the approximated density $p^r(r)/Z_r$ given in \cite[eqn (33)]{DNR}. 
\end{remark}

\begin{remark}  Theorem \ref{thm sigma0} deals with small noise when $\mu > 0$.  For completeness we give small noise results for $\mu < 0$ and $\mu = 0$ also.  If $\mu < 0$ it follows from \eqref{int r2} that $\E[r^2] \sim \sigma^2/(-\mu)$ as $\sigma \to 0$ and so $\lambda(\mu,a,b,\sigma) = \mu  + {\cal O}(\sigma^2)$ as $\sigma \to 0$.
If $\mu = 0$ then scaling \eqref{scale2} implies $\lambda(0,a,b,\sigma) = \sqrt{a}\sigma \lambda(0,1,b/a,1)$.  Proposition \ref{prop Jhat} shows that $\lambda(0,1,b/a,1) < 0$ for sufficiently small $b/a$, and Theorem \ref{thm binf} shows $\lambda(0,1,b/a,1)> 0$ for sufficiently large $b/a$.
\end{remark}

\begin{corollary}  \label{cor sigma0} (i) Fix $\mu > 0$ and $a > 0$.  For all $k >0$ there exists $\sigma_0 > 0$ such that $\lambda(\mu,a,b,\sigma) < 0$ whenever $0 \le b \le k$ and $0 < \sigma \le \sigma_0$.

(ii) Fix $a >0$ and $\sigma >0$.  For all $k >0$ there exists $\mu_0 > 0$ such that $\lambda(\mu,a,b,\sigma) < 0$ whenever $0 \le b \le k$ and $\mu \ge \mu_0$.
   \end{corollary} 

The second part of Corollary \ref{cor sigma0} uses \eqref{scale3} in Proposition \ref{prop scale} to convert a small $\sigma$ result into a large $\mu$ result.  Combining this with Proposition \ref{prop Jhat} gives

\begin{corollary} \label{cor sigma0 2}  %(i) 
There exists $k_0 > 0$ such that $\lambda(\mu,a,b,\sigma) <0$ whenever $0 \le b/a  \le k_0$. 

%(ii) For all $a > 0$, $b \ge 0$ and $\sigma > 0$ there exist $\mu_1 < \mu_2$ such that $\lambda(\mu,a,b,\sigma) < 0$ when $\mu \le \mu_1$ and and $\mu \ge \mu_2$.
\end{corollary}

\n{\bf Proof.% of Corollary \ref{cor sigma0 2}.
} %(i) 
By \eqref{scale2} in Proposition \ref{prop scale} it suffices to consider the case $a=\sigma=1$.  Choose $b_0 > 0$, and then let $\mu_0$ be as in Corollary \ref{cor sigma0} with $k = b_0$.  Define $b_1= \widehat{J}(\mu_0/\sqrt{2})$ and let $k_0 = \min(b_0,b_1)$.  Now suppose $b \le k_0$.  If $\mu \le \mu_0$ then $b \le b_ 1 = \widehat{J}(\mu_0/\sqrt{2}) \le \widehat{J}(\mu/\sqrt{2})$ so that $\lambda(\mu,1,b,1) < 0$ by Proposition \ref{prop Jhat}, and if $\mu \ge \mu_0$ then $b \le b_0$ and so $\lambda(\mu,1,b,1) <0$ by Corollary \ref{cor sigma0}. 
%(ii) The existence of $\mu_1$ with the required property follows from Proposition \ref{prop Jhat} and the existence of $\mu_2$ with the required property follows from Corollary \ref{cor sigma0}.
 \qed

\subsection{Stochastic bifurcation} \label{sec stocbif} When $\sigma = 0$ the system reduces to the deterministic Hopf bifurcation with a qualitative change in the dynamics as the bifurcation parameter $\mu$ passes through 0.  Suppose now $a > 0$, $b \ge 0$ and $\sigma > 0$ are fixed, and consider the possible changes in stability as $\mu$ is varied.  We have seen
   \begin{description}
      
   \item[\quad (i)]   If $0 \le b/a \le k_0$ then $\lambda <0$ for all $\mu$.  (This is Corollary \ref{cor sigma0 2}.)
   
   \item[\quad (ii)]   There exist $\mu_1 < \mu_2$ such that $\lambda < 0$ when $\mu \le \mu_1$ and when $\mu \ge \mu_2$.  (The result for $\mu \le \mu_1$ follows from Proposition \ref{prop Jhat} and the result for $\mu \ge \mu_2$ follows from Corollary \ref{cor sigma0}.)
   
   \item[\quad (iii)]  If $b/a \ge k_1$ then $\lambda > 0$ for $|\mu| \le c \sqrt{a} \sigma$.  (This is Corollary \ref{cor binf unif2}.)
   \end{description} 

Note that the constants $k_0$, $k_1$ and $c$ do not depend on any of the parameter values.  The property {\bf (i)} implies that when the twist factor $b/a$ is small then additive noise destroys the Hopf bfurcation.   But when the twist factor is sufficiently large then {\bf (ii)} and {\bf (iii)} together imply that $\lambda$ changes sign from negative to positive and back to negative as $\mu$ increases from $-\infty$ through 0 to $\infty$.

\subsection{Stability diagram} \label{sec stab diag}

Recall from \eqref{b} that $\lambda(\mu,a,b,\sigma) = \lambda(\mu,a,-b,\sigma)$, and from Proposition \ref{prop scale} that 
   $$ \lambda(\mu,a,b,\sigma) = \sigma\sqrt{a} \,\lambda\left(\frac{\mu}{\sigma \sqrt{a}},1,\frac{b}{a},1\right).
   $$ 
Therefore for the stability diagram it is enough to consider the sign of $\lambda(\mu,1,b,1)$ for $(\mu,b)$ in the half space $\R \times [0,\infty)$.

  \begin{figure}[h] 
\begin{center}
 \includegraphics[scale=0.9]{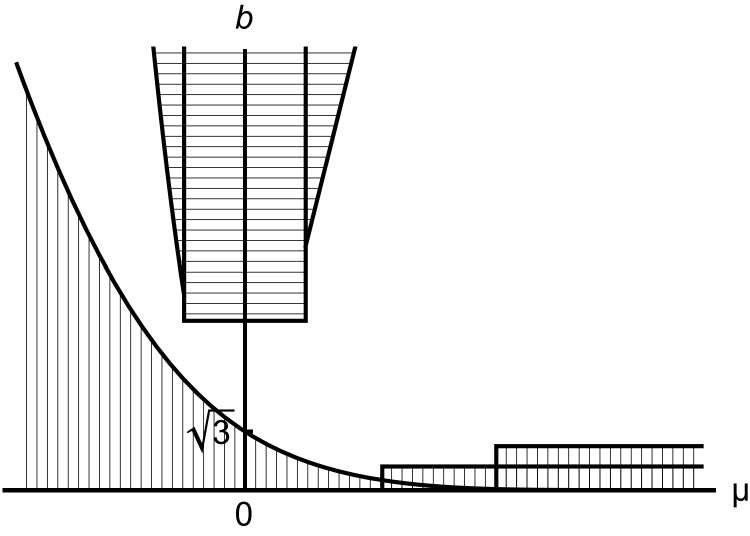}
 \end{center}
 \caption{Stability diagram: horizontally shaded regions where $\lambda(\mu,1,b,1) > 0$, and vertically shaded regions where $\lambda(\mu,1,b,1) < 0$}\label{fig stab diag}
 \end{figure}

Figure \ref{fig stab diag} shows horizontally shaded regions where $\lambda(\mu,1,b,1) > 0$, as given by Corollary \ref{cor binf unif2}.  It also shows vertically shaded regions where $\lambda(\mu,1,b,1) < 0$, given by the set $b \le \widehat{J}(\mu/\sqrt{2})$ in Proposition \ref{prop Jhat} together with two applications of Corollary \ref{cor sigma0}(ii) with different values for $k$.  We note that the vertically shaded regions in Figure \ref{fig stab diag} provide a ``visual proof'' of Corollary \ref{cor sigma0 2}. 

\s

  For large negative $\mu$ the gap between the stable and unstable regions shown is bounded by curves $b = \widehat{J}(\mu/\sqrt{2})$ and $b = k_3\mu^2$ for some $k_3$.  Since $\widehat{J}(\mu/\sqrt{2}) \sim \mu^2$ as $\mu \to -\infty$, both curves have asymptotically quadratic growth, but the coefficients 1 and $k_3$ may be very far apart.   
  For large positive $\mu$ the situation is rather different.  In Corollary \ref{cor sigma0}(ii) there is no simple description of the dependence of $\mu_0$ on $k$, and hence no simple description of the horizontal displacements of the rectangular blocks on the lower right as their heights are increased.   Instead we have the following strengthening of Chemnitz and Engel's Theorem \ref{thm CE}.

\begin{theorem} \label{thm CE unif}    Fix $\mu > 0$ and $a >0$ and a compact set $D \subset (0,\infty)$.  For all $\e > 0$ there exists $\sigma_0 > 0$ such that 
    \begin{equation} \label{CE unif}
 \left| \lambda(\mu,a,b,\sigma) - 2 \mu \Psi\left(\frac{b^2 \sigma^2 }{2 \mu^2 a}\right) \right| <\e
     \end{equation}
whenever $b\sigma \in D$ and $0 < \sigma \le \sigma_0$ 
\end{theorem}   

The proof of Theorem \ref{thm CE unif} is given in Section \ref{sec CE proof}.   Using \eqref{scale3}, and recalling that $\Psi(\zeta) < 0$ for $0 < \zeta < c_\ast$ and $\Psi(\zeta) > 0$ for $\zeta > c_\ast$, this gives     

\begin{corollary} \label{cor CE unif} (i) Given $0 < c_1 < c_2 < \sqrt{2c_\ast}$, there exists $\mu_0 > 0$ such that $\lambda(\mu,1,b,1) <0$ whenever $\mu \ge \mu_0$ and $c_1 \le b/\mu \le c_2$.

(ii) Given $\sqrt{2c_\ast} < c_3 < c_4$, there exists $\mu_0 > 0$ such that $\lambda(\mu,1,b,1) >0$ whenever $\mu \ge \mu_0$ and $c_3 \le b/\mu \le c_3$. 
   \end{corollary}

  \begin{figure}[t] 
\begin{center}
 \includegraphics[scale=0.9]{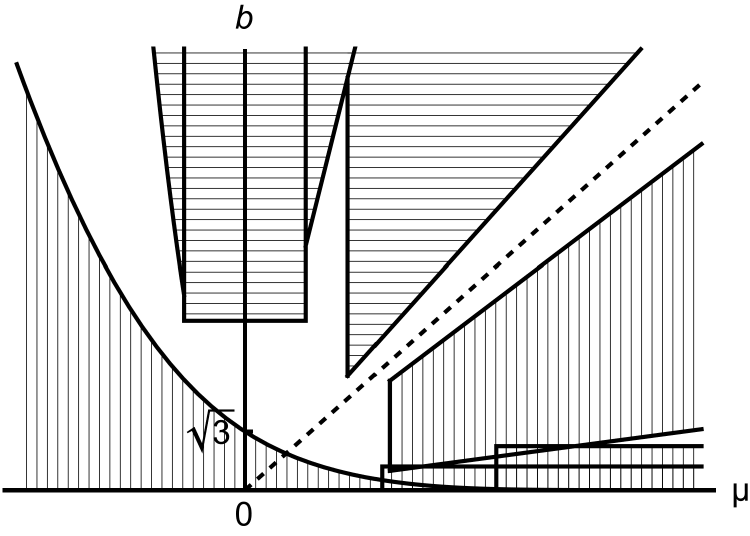}
 \end{center} 
 \caption{Stability diagram: horizontally shaded regions where $\lambda(\mu,1,b,1) > 0$, and vertically shaded regions where $\lambda(\mu,1,b,1) < 0$. Dashed line is $b/\mu = \sqrt{2c_\ast} \approx 2.66$.}\label{fig stab diag2}
 \end{figure} 
Figure \ref{fig stab diag2} shows the effect of Corollary \ref{cor CE unif} on the stability diagram shown in Figure \ref{fig stab diag}.  The dashed line is $b/\mu = \sqrt{2c_\ast} \approx 2.66$.  Note that $c_4$ can be chosen to be the $k_2$ from Corollary \ref{cor binf unif} so there is no gap in the unstable region for large $b/\mu$.  However, without a further strengthening of Corollary \ref{cor sigma0}, there may be gaps in the stable region for small $b/\mu$.  

\subsubsection{Numerical stability diagram} \label{sec stab diag num}

\begin{figure}[h]
\begin{center} \includegraphics[scale=0.9]{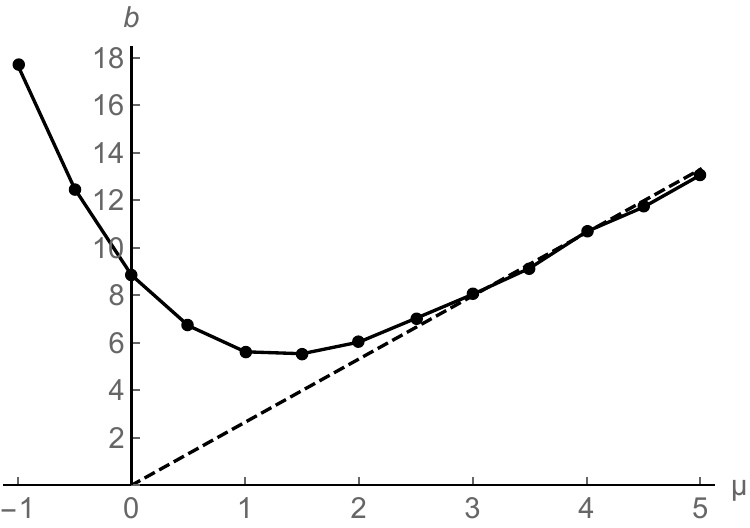} \end{center}
\caption{Numerical approximation of $\lambda(\mu,1,b,1) = 0$.  Dashed line is $b/\mu = \sqrt{2c_\ast} \approx 2.66$.}\label{fig num}
 \end{figure}

Figures \ref{fig stab diag} and \ref{fig stab diag2} are symbolic representations of theoretical results.  In contrast Figure \ref{fig num} is the result of numerical simulation.  Instead of using the Furstenburg-Khasminskii formula \eqref{lam khas}, which involves the time costly evaluation of trigonometric functions, the original system (\ref{X},\ref{U}) was simulated using the Euler method with time step $h=0.0001$.  Then the finite-time Lyapunov exponent (see for example \cite{DELR}) at time $n=10,000$ was used as a proxy for the true value of the Lyapunov exponent.   In Figure \ref{fig num} the unstable region $\lambda(\mu,1,b,1) > 0$ in $(\mu,b)$ parameter space lies above the curve and the stable region $\lambda(\mu,1,b,1) < 0$ lies below the curve.  The dashed line is $b/\mu =\sqrt{2c_\ast} \approx 2.66$.  Recall from Corollary \ref{cor CE unif} that it gives the asymptotic line of separation between the regions as $\mu \to \infty$. It is notable that large $\mu$ theory agrees with the numerics for relatively small values of $\mu$.

There are several previous papers with similar numerical computations of $\lambda$.  
%Figures 8 and 9 of Wieczorek \cite{Wie} show the curves $\lambda(\Lambda,1,\alpha,\sqrt{D_{ext}}) = 0$ for positive $\Lambda$.  
The straight lines for $\lambda= 0$ in Wieczorek \cite[Figure 8]{Wie} and the threefold repetition of the same curve in  of \cite[Figure 9]{Wie} are direct consequences of the scaling \eqref{scale2}.  After converting the notation in \cite[Figure 9]{Wie} we see for example $\lambda(\mu,1,6,1) = 0$ for $\mu \approx 2.0$ and $0.79$.  Also the horizontal asymptotic $\alpha \approx 9$ in \cite[Figure 9]{Wie} converts to $\lambda(0,1,b,1) = 0$ for $b \approx 9$.
DeVille, Namachchivaya and Rapti \cite[Figure 1]{DNR} show $\lambda(\mu,1,b,1)$ for the region 
$1 \le \mu \le 8$ and $1 \le b \le 20$, and the plot shows for example $\lambda(2,1,b,1) = 0$ for $b \approx 6.4$ and $\lambda(4,1,b,1) = 0$ for $b \approx 12$. 
Doan, Engel, Lamb and Rasmussen \cite[Figure 3]{DELR} show $\lambda(\mu,1,b,1)$ for the region $-0.2 \le \mu \le 0.2$ and $8 \le b \le 10.5$ in which the curve $\lambda(\mu,1,b,1) = 0$ is close to linear, and the plot shows $\lambda(0,1,b,1) = 0$ for $b \approx 9.1$.
  
%These numerical observations are consistent with the results in our bifurcation diagram above.   

\section{Furstenberg-Khasminskii formula} \label{sec FK}

We take advantage of the rotational symmetry of the system \eqref{X} by using polar coordinates $X_t = r_t \begin{bmatrix}\cos \phi_t \\ \sin \phi_t \end{bmatrix}$.  For matrix calculations it is convenient to introduce the infinitesimal rotation matrix $J = \begin{bmatrix} 0 &-1  \\ 1 & 0 \end{bmatrix}$.  Following \cite{DNR, CE23}, define new independent Wiener processes
  \begin{align*}
     W_t^r &  =  \int_0^t \cos \phi_s dW_s^1 + \int_0^t \sin \phi_s dW_s^2, \\
     W_t^\phi & =  -\int_0^t \sin \phi_s dW_s^1 + \int_0^t \cos \phi_s dW_s^2.
     \end{align*} 
It\^o's formula applied to \eqref{X} gives 
   \begin{align} \label{r bis}
   dr_t & = \left(\mu r_t- ar_t^3+ \frac{\sigma^2}{2r_t}\right)dt+ \sigma dW_t^r,\\
   d \phi_t &  =  (\omega+br_t^2)dt + \frac{\sigma}{r_t}dW_t^\phi. \nonumber
   \end{align}
Putting $X_t = r_t \begin{bmatrix}\cos \phi_t \\ \sin \phi_t \end{bmatrix}$ in \eqref{U} gives
   \begin{equation} \label{U2}
   dU_t = \left( \begin{bmatrix} \mu  & -\omega\\
  \omega & \mu \end{bmatrix} + r_t^2\begin{bmatrix} -a & -b \\ b & -a \end{bmatrix}\begin{bmatrix} 2+\cos 2 \phi_t & \sin 2 \phi_t \\ \sin 2 \phi_t& 2-\cos 2 \phi_t \end{bmatrix}\right)U_t dt
  \end{equation}
Write $U_t = \|U_t\| \begin{bmatrix} \cos \theta_t\\ \sin \theta_t\end{bmatrix}$.  Then It\^o's formula applied to \eqref{U2} gives  
  \begin{align}
  d \log \|U_t\|  &  =  \frac{1}{\|U_t\|^2} \langle U_t,dU_t \rangle \nonumber \\
  & = \left\langle  \begin{bmatrix} \cos \theta_t \\ \sin \theta_t \end{bmatrix}, \left( \begin{bmatrix} \mu  & -\omega\\
  \omega & \mu \end{bmatrix} + r_t^2\begin{bmatrix} -a & -b \\ b & -a \end{bmatrix}\begin{bmatrix} 2+\cos 2 \phi_t & \sin 2 \phi_t \\ \sin 2 \phi_t& 2-\cos 2 \phi_t \end{bmatrix}\right) \begin{bmatrix} \cos \theta_t \\\sin \theta_t \end{bmatrix} \right\rangle dt \nonumber \\
  & := \widehat{Q}(r_t,\phi_t,\theta_t)dt \label{Qhat}
  \end{align}
and similarly 
  \begin{align}
  d \theta_t &  =  \frac{1}{\|U_t\|^2}\langle JU_t,dU_t \rangle \nonumber \\
  & = \left\langle  \begin{bmatrix} -\sin \theta_t \\ \cos \theta_t \end{bmatrix}, \left( \begin{bmatrix} \mu  & -\omega\\
  \omega & \mu \end{bmatrix} + r_t^2\begin{bmatrix} -a & -b \\ b & -a \end{bmatrix}\begin{bmatrix} 2+\cos 2 \phi_t & \sin 2 \phi_t \\ \sin 2 \phi_t& 2-\cos 2 \phi_t \end{bmatrix}\right) \begin{bmatrix} \cos \theta_t \\\sin \theta_t \end{bmatrix} \right\rangle dt \nonumber \\
  & := \widehat{H}(r_t,\phi_t,\theta_t) dt, \label{Hhat}
 \end{align}
say.  %for certain functions $\widehat{Q}$ and $\widehat{F}$. 
Now define $\psi = \theta-\phi$.  A direct calculation gives  
    \begin{align} \nonumber
    \widehat{Q}(r,\theta,\phi) & =  \mu-2ar^2 +r^2\Bigl( b \sin (2\theta-2\phi)- a \cos (2\theta-2\phi) \Bigr)\\
    & = \mu-2ar^2 +r^2\Bigl( b \sin 2\psi- a \cos 2 \psi \Bigr):=Q(r,\psi) \label{QQ}
    \end{align}
and   \begin{align*} 
   \widehat{H}(r,\theta,\phi) & = \omega+2br^2  +r^2\Bigl( a \sin (2\theta-2\phi)+ b \cos (2\theta-2\phi) \Bigr)\\
    & = \omega+2br^2  +r^2\Bigl( a \sin 2\psi + b \cos 2\psi \Bigr):=H(r,\psi).  
   \end{align*}  
Then   
  \begin{align} \nonumber
   d \psi_t   =  d\theta_t-d\phi_t & =  H(r_t,\psi_t)dt -\left(\omega+br_t^2 \right)dt - \frac{\sigma}{r_t}dW_t^\phi \nonumber \\
    & =  2r^2 \cos \psi (b \cos \psi+a \sin \psi)dt- \frac{\sigma}{r_t}dW_t^\phi . \label{psi bis}
   \end{align}
The equations \eqref{r bis} for $r_t$ and \eqref{psi bis} for $\psi_t$ together give the system (\ref{r},\ref{psi}) with generator $\cal L$ described in Section \ref{sec eval}, and \eqref{QQ} gives the formula \eqref{Q} for $Q(r,\psi)$.  Integrating \eqref{Qhat} gives 
    \begin{equation} \label{norm}
   \frac{1}{t} \log\|U_t\| = \frac{1}{t}\log \|U_0\| + \frac{1}{t} \int_0^t Q(r_s,\psi_s)ds.
  \end{equation}
Since the operator ${\cal L}$ is elliptic with invariant probability $\widetilde{\nu}$ on $M$, the Furstenberg-Khasminskii formula \eqref{lam khas} follows by the almost-sure ergodic theorem.  Note that the limit on the right exists almost surely and takes the same value for all starting points $X_0$ and all non-zero initial vectors $U_0$. 
  
\subsection{Adjoint method} \label{sec adj}

 The essential idea in the adjoint method, introduced by Arnold, Papanicolaou and Wihstutz \cite{APW86}, is to replace the evaluation of the integral $\int_M Q(r,\psi)d\widetilde{\nu}(r,\psi)$ in the Furstenberg-Khasminskii formula by the search for a function $F(r,\psi)$ and a constant $\lambda$ satisfying ${\cal L}F(r,\psi) = Q(r,\psi)-\lambda$, or equivalently $Q(r,\psi)-{\cal L}F(r,\psi) = \lambda$.  Under suitable integrability and growth conditions, see Proposition \ref{prop-bg}, we have $\int_M {\cal L}F(r,\psi)d\widetilde{\nu}(r,\psi) = 0$ and so 
  $$
  \int_M Q(r,\psi)d\widetilde{\nu}(r,\psi) = \int_M \bigl(Q(r,\psi)-{\cal L}F(r,\psi)\bigr)d\widetilde{\nu}(r,\psi) = \int_M \lambda\, d\widetilde{\nu}(r,\psi) = \lambda.
  $$
When applying the method it is often impossible to find a function $F$ for which $Q(r,\psi)-{\cal L}F(r,\psi) $ is constant, but it is frequently possible to find a function $F$ for which $Q(r,\psi)-{\cal L}F(r,\psi) $ is easier to analyze than the original $Q(r,\psi)$.  The proof of Theorem \ref{thm binf} uses the adjoint method twice.  The first step uses a function $F$ motivated by the heuristic argument in Section \ref{sec h binf} to exhibit the $b^{2/3}$ growth rate, see Proposition \ref{prop QLF}, and the second step uses a function $G$ to remove the $\psi$ dependence from the problem, see Lemma \ref{lem G} and Proposition \ref{prop Phi}.  The proof of Theorem \ref{thm CE unif} uses a very similar two step method, starting with exactly the same $F$ as in the proof of Theorem \ref{thm binf}.  In contrast, the proof of Theorem \ref{thm sigma0} starts in a more conventional manner, looking for an asymptotic expansion $F(r,\psi) = F_0(r,\psi)+\sigma^2 F_2(r,\psi) +\cdots$ in powers of the small parameter $\sigma^2$.  The novelty here is in the construction of the functions $F_0$ and $F_2$.  Unlike the original setting of \cite{APW86} involving an elliptic operator on a compact state space, the equations (\ref{adj},\ref{adj2}) for $F_0$ and $F_2$ involve a first order differential operator on the non-compact space $M = (0,\infty) \times \R/(\pi \Z)$.  The functions $F_0$ and $F_2$ have singularities at the hyperbolic fixed point of the corresponding vector field, and the proof requires an {\it a priori} estimate on the invariant probability measure $\widetilde{\nu}$ near this fixed point as $\sigma \to 0$, see Lemma \ref{lem nu U}.

\section{Heuristic arguments}  \label{sec h}
Before presenting the rigorous proofs in Section \ref{sec binf proof} onwards, we give here some heuristic arguments for the results in Theorems \ref{thm binf} - \ref{thm CE unif}.  These are seen most easily when the derivative process $U_t$ is written in terms of the rotating %orthonormal 
frame  given by the unit radial and tangential vectors $\begin{bmatrix} \cos \phi_t \\\sin \phi_t \end{bmatrix}$ and $\begin{bmatrix} -\sin \phi_t  \\ \cos \phi_t \end{bmatrix}$ at the point $X_t$.  
For $\tau \in \R$ define $R_\tau = \begin{bmatrix}\cos \tau & -\sin \tau \\ \sin \tau & \cos \tau \end{bmatrix}$, corresponding to rotation through the angle $\tau$. Note that $R'_\tau = JR_\tau$ where $J = \begin{bmatrix} 0 & -1 \\ 1 & 0 \end{bmatrix}$. 
The components of $U_t$ relative to the rotating frame are given by the vector $V_t$, say, where 
             $
             U_t = R_{\phi_t} V_t.
             $
So $V_t = R_{-\phi_t}U_t$ and then 
   \begin{align}
      dV_t & = R_{-\phi_t}dU_t-JR_{-\phi_t} U_t d \phi_t  -\frac{1}{2}R_{-\phi_t}U_t(d\phi_t)^2 \nonumber  \\
        & = R_{-\phi_t} \left( \begin{bmatrix} \mu  & -\omega\\
          \omega & \mu \end{bmatrix} + r_t^2\begin{bmatrix} -a & -b \\ b & -a \end{bmatrix}
            \begin{bmatrix} 2+ \cos 2 \phi_t & \sin 2 \phi_t \\ \sin 2 \phi_t& 2-\cos 2 \phi_t \end{bmatrix}      \right)R_{\phi_t}V_t dt \nonumber \\
            & \quad -(\omega+br_t^2)J V_t dt - \frac{\sigma}{r_t}JV_tdW_t^\phi -\frac{\sigma^2}{2r_t^2}V_t dt.  \label{dV}
   \end{align}              
Since 
 $$R_{-\phi_t}\begin{bmatrix} \cos 2 \phi_t & \sin 2 \phi_t\\ \sin 2 \phi_t & -\cos 2 \phi_t \end{bmatrix}R_{\phi_t} = \begin{bmatrix} 1 & 0 \\ 0 & -1 \end{bmatrix}  \mbox{ and }  
R_{-\phi_t}J \begin{bmatrix} \cos 2 \phi_t & \sin 2 \phi_t\\ \sin 2 \phi_t & -\cos 2 \phi_t \end{bmatrix}R_{\phi_t}
=  \begin{bmatrix} 0 &1  \\1 & 0  \end{bmatrix}
$$ 
the right side of \eqref{dV} simplifies and we get
   \begin{align}
      dV_t & = \left( \begin{bmatrix} \mu  & -\omega\\
          \omega & \mu \end{bmatrix} + 2r_t^2\begin{bmatrix} -a & -b \\ b & -a \end{bmatrix}
            -a r_t^2 \begin{bmatrix} 1 & 0 \\0 & -1 \end{bmatrix} + br_t^2 \begin{bmatrix} 0 & 1 \\ 1 & 0 \end{bmatrix} \right)V_t dt  \nonumber \\
            & \quad -(\omega+br_t^2)\begin{bmatrix} 0 & -1 \\ 1 & 0 \end{bmatrix} V_t dt - \frac{\sigma}{r_t}\begin{bmatrix} 0 & -1 \\ 1 & 0 \end{bmatrix}V_tdW_t^\phi  -\frac{\sigma^2}{2r_t^2}V_t dt\nonumber \\
            & =  \begin{bmatrix} \mu - 3ar_t^2 - \sigma^2/2r_t^2  & 0 \\ 2br_t^2 & \mu - ar_t^2-\sigma^2/2r_t^2 \end{bmatrix} V_t dt + \frac{\sigma}{r_t} \begin{bmatrix} 0 & 1 \\ -1 & 0 \end{bmatrix} V_t  dW_t^\phi.   \label{dV2}
   \end{align} 
This equation appears in Stratonovich form as \cite[eqn (2.12)]{CE23}. 

\begin{remark}Since $V_t$ is obtained from $U_t$ by rotation through the angle $-\phi_t$, it follows that $V_t = \|U_t\| \begin{bmatrix} \cos \psi_t\\ \sin \psi_t\end{bmatrix}$, where $\psi_t = \theta_t-\phi_t$ is defined in Section \ref{sec FK}.   Applying the methods of equations (\ref{Qhat},\ref{Hhat}) to the It\^{o} SDE for $V_t$ in place of the ordinary differential equation for $U_t$ gives an alternative method for deriving the function $Q(r,\psi)$ and the equation \eqref{psi} for $\{\psi_t: t \ge 0\}$.   The calculation is a direct application of It\^o's formula, and is similar to those carried out for Stratonovich versions in \cite[eqns (5,6)]{IL01} and \cite[Lemma 3.2]{CE23}.

\end{remark}

\subsection{Heuristic argument for Theorem \ref{thm binf}}  \label{sec h binf} Fix $\mu \in\R$, $a >0$ and $\sigma >0$.  As $b \to \infty$ we see the emergence of a nilpotent structure in the drift matrix of \eqref{dV2}, and this suggests the transformation used by Pinsky and Wihstutz \cite{PW88} for nilpotent linear SDEs.  In a linear SDE with constant coefficients we would consider $\widetilde{V}_t = \begin{bmatrix}A & 0 \\ 0 & 1 \end{bmatrix} V_t$ where $A = Bb^{1/3}$ for a suitable constant $B$. This gives
\begin{align}
     d\widetilde{V}_t
       & =  \nonumber  
       \begin{bmatrix} \mu - 3ar_t^2 - \sigma^2/2r_t^2  & 0 \\  2b^{2/3}r_t^2/B & \mu - ar_t^2-\sigma^2/2r_t^2 \end{bmatrix}\widetilde{V}_t dt \\[1ex]
       & \hspace{25ex}  
       + \frac{\sigma}{r_t}\begin{bmatrix} 0 & Bb^{1/3} \\ -1/(Bb^{1/3}) & 0 \end{bmatrix}\widetilde{V}_t dW_t^\phi \label{AUV} \\[2ex]
         & = \nonumber b^{2/3} \begin{bmatrix} 0 & 0 \\ 2r_t^2/B & 0 \end{bmatrix}\widetilde{V}_tdt +b^{1/3}\begin{bmatrix} 0 & B  \sigma /r_t \\ 0 & 0 \end{bmatrix}\widetilde{V}_tdW_t^\phi   \nonumber \\[1ex]
         &  \qquad  + \mbox{ terms of lower order in }b \mbox{ as }b \to \infty.  \label{AUVI}
              \end{align} 
Recall
  $$
   dr_t  = \bigl(\mu r_t - ar_t^3+\frac{\sigma^2}{2r_t}\bigr)dt + \sigma dW_t^r.
   $$
   We see a separation of time scales: $\widetilde{V}_t$ moves at a rate proportional to $b^{2/3}$ whereas $r_t$ moves at rate independent of $b$.   Proceeding under the heuristic that the fast motion in \eqref{AUVI} can be analysed while treating $r_t$ as a constant, we choose $B$ so that $2r^2/B = (B  \sigma /r)^2$.  This gives $B = (2r^4/\sigma^2)^{1/3}$ and then two leading terms in \eqref{AUVI} correspond to a time changed version of the equation \eqref{YA}:
  $$
 dY_t = \begin{bmatrix} 0 & 0 \\1 & 0 \end{bmatrix}Y_tdt
   +\begin{bmatrix} 0 & 1 \\0 & 0 \end{bmatrix}Y_tdW_t
  $$
with time change factor $b^{2/3}\times(2r^2/B) = (2 b \sigma  r)^{2/3}$.  The time change gives a Lyapunov exponent $(2 b \sigma  r)^{2/3}\gamma_0$.  This calculation was carried out under the assumption that $r$ is fixed.  The final stage in this heuristic is to average over $r$ with respect to the invariant measure $\nu_{(\mu,a,\sigma)}$.  Thus we would get
    $$
    \lambda(\mu,a,b,\sigma) \sim (2 b \sigma)^{2/3}\gamma_0 \int_{\R_+} r^{2/3}d\nu_{(\mu,a,\sigma)}(r)
    $$
as $b \to \infty$.  However this is non-rigorous.  The fact that $B$ is a function of $r=r_t$ introduces extra terms into \eqref{AUV} and hence extra lower order terms in \eqref{AUVI}.  In the setting of the Furstenberg-Khasminskii formula, there are lower order terms arising from the $Q$ function, and there are lower order terms in the generator ${\cal L}$.  All of these remainder terms need to be controlled, but the analysis of the remainder terms involves functions of $r$ which are not integrable near 0.  Instead of pursuing this idea directly, we use it as motivation for the proof of Theorem \ref{thm binf} in Section \ref{sec binf proof}.

More precisely, if we make the change from $V$ to $\widetilde{V} = \begin{bmatrix} A(r) & 0 \\0 & 1 \end{bmatrix}V$ with $A(r) = (2 b r^4/\sigma^2)^{1/3}$ and then apply the method of Khasminskii to the SDE for $\widetilde{V}_t$ then the integrand in the Furstenberg-Khasminskii formula changes from
$Q(r,\psi)$ to 
    $$
    Q(r,\psi) +{\cal L} \log\left(\frac{\|\widetilde{V}\|}{\|V\|}\right) = 
    Q(r,\psi) + \frac{1}{2}{\cal L}\Bigl(\log\bigl(A(r)^2 \cos^2\psi+\sin^2 \psi\bigr)\Bigr).
    $$
However the function $(r,\psi) \mapsto \log \bigl(A(r)^2 \cos^2 \psi + \sin^2 \psi\bigr)$ blows up as $r \to 0$ and $\psi \to 0$.  In Section \ref{sec binf proof} we keep the original Furstenberg-Khasminskii formula \eqref{lam khas}, and then use a variation of the adjoint method by replacing $Q(r,\psi)$ with $Q(r,\psi) - {\cal L}F(r,\psi)$ where the function 
    $$
    F(r,\psi) = -\frac{1}{2} \log\bigl(A(r)^2 \cos^2 \psi + 1\bigr) = \frac{1}{2} \log \left(\frac{1+\tan^2 \psi}{A(r)^2+1+\tan^2 \psi}\right)
    $$
is better behaved as $r \to 0$.  The change of variable from $z = \tan \psi$ to $w = z/A(r)$ used in Section \ref{sec binf proof} corresponds exactly to the change from $V$ to $\widetilde{V}=  \begin{bmatrix}A(r) & 0 \\ 0 & 1 \end{bmatrix} V$.

\subsection{Heuristic argument for Theorem \ref{thm CE unif}} \label{sec h CE} When $b\sigma$ is bounded above and below and $\sigma \to 0$ the hierarchy of terms in \eqref{AUV} is changed.  Taking $B = (2r^4/\sigma^2)^{1/3}$ again we have
  \begin{align*}
     d\widetilde{V}_t
       & =   \begin{bmatrix} \mu - 3ar_t^2   & 0 \\  (2b \sigma r_t)^{2/3} & \mu - ar_t^2 \end{bmatrix}\widetilde{V}_t dt+ \begin{bmatrix} 0 & (2b\sigma r_t)^{1/3} \\ 0 & 0 \end{bmatrix} \widetilde{V}_t dW_t^\phi\\[1ex]
                 &  \qquad  + \mbox{ terms of lower order in }\sigma. %\mbox{ as }\sigma \to 0. 
              \end{align*} 
Moreover since $r_t \to \sqrt{\mu/a}$ in distribution as $\sigma \to 0$ we get
    \begin{align}
     d \widetilde{V}_t
       & = \nonumber  \begin{bmatrix} -2\mu   & 0 \\  (2b \sigma \sqrt{\mu/a})^{2/3} & 0 \end{bmatrix}\widetilde{V}_t dt+ \begin{bmatrix} 0 & (2b\sigma \sqrt{\mu}{a})^{1/3} \\ 0 & 0 \end{bmatrix} \widetilde{V}_t dW_t^\phi\\[1ex]
       &  \nonumber \qquad +\mbox{ terms with a factor }(r_t-\sqrt{\mu/a}) \,\, + \mbox{ terms of lower order in }\sigma\\[1ex]
        & = \nonumber  2 \mu \begin{bmatrix} -1   & 0 \\  (b^2\sigma^2/2\mu^2 a)^{1/3} & 0  \end{bmatrix}\widetilde{V}_t dt+ \sqrt{2\mu}\begin{bmatrix} 0 & (b^2\sigma^2/2\mu^2 a)^{1/6} \\ 0 & 0 \end{bmatrix} \widetilde{V}_t dW_t^\phi\\[1ex]
                 &  \qquad   +\mbox{ terms with a factor }(r_t-\sqrt{\mu/a}) \,\, + \mbox{ terms of lower order in }\sigma%+ \mbox{ smaller terms as }\sigma \to 0
                 .  \label{AUVII}
              \end{align} 
The two leading terms in \eqref{AUVII} correspond to a time changed version of the equation \eqref{YB}:
  $$
  dY_t = \begin{bmatrix} -1  & 0 \\ \zeta^{1/3} & 0 \end{bmatrix}Y_t dt + \begin{bmatrix} 0 & \zeta^{1/6} \\ 0 & 0 \end{bmatrix} Y_t dW_t 
     $$ with $\zeta = b^2\sigma^2/2\mu^2 a$ and time change factor $2 \mu$.   Since \eqref{YB} has Lyapunov exponent $\Psi(\zeta) = \Psi\left(\dfrac{b^2\sigma^2}{2 \mu^2 a}\right)$, then the time change gives a Lyapunov exponent $2\mu \Psi\left(\dfrac{b^2\sigma^2}{2 \mu^2 a}\right)$.  Again this is a non-rigorous argument.  Sionce the heuristic argument starts with the same transformation $V \mapsto \widetilde{V}$ as in Section \ref{sec h binf}, the proof of Theorem \ref{thm CE unif}, given in Section \ref{sec CE proof}, starts with the same adjoint method and uses the same Corollary \ref{cor lam} that is used in the proof of Theorem \ref{thm binf}. 

\subsection{Heuristic argument for Theorem \ref{thm sigma0}} \label{sec h sigma0}

%In the rotating frame the linearized process is given by the SDE
%    \begin{equation}
%      \begin{bmatrix}d U_t \\dV_t \end{bmatrix}
%          =  \begin{bmatrix} \mu - 3ar_t^2 - \sigma^2/2r_t^2  & 0 \\ 2br_t^2 & \mu - ar_t^2-\sigma^2/2r_t^2 \end{bmatrix}\begin{bmatrix} U_t \\V_t \end{bmatrix}dt + \frac{\sigma}{r_t}\begin{bmatrix} 0 & 1 \\ -1 & 0 \end{bmatrix} \begin{bmatrix} U_t \\V_t \end{bmatrix}dW_t^4.  \label{UV}
%              \end{equation} 
As $\sigma \to 0$ then $r_t \to \sqrt{\mu/a}$ in distribution.  Replacing $r_t$ with $\sqrt{\mu/a}$ in the SDE \eqref{dV2} for $V_t$ gives the constant coefficient SDE
  \begin{equation}
     dV_t
          =  \begin{bmatrix} -2\mu - \sigma^2a/2\mu  & 0 \\ 2b\mu/a & - \sigma^2 a/2\mu \end{bmatrix}V_t dt + \frac{\sigma \sqrt{a}}{\sqrt{\mu}} \begin{bmatrix} 0 & 1 \\ -1 & 0 \end{bmatrix} V_t dW_t^\phi.  \label{dV3}
              \end{equation} 
The top Lyapunov exponent for \eqref{dV3} is $-\sigma^2 a/2\mu$ plus the top Lyapunov exponent of the SDE
   \begin{equation}
     dV_t
          =  \begin{bmatrix} -2\mu   & 0 \\ 2b\mu/a & 0 \end{bmatrix} V_t dt + \frac{\sqrt{a}\sigma}{\sqrt{\mu}} \begin{bmatrix} 0 & 1 \\ -1 & 0 \end{bmatrix} V_t dW_t^\phi.  \label{dV4}
              \end{equation} 
%The asymptotic as $\sigma \to 0$ for the top Lyapunov exponent of \eqref{UV3} is given by .  
Changing basis $ \widehat{V}_t  = \begin{bmatrix}1/a & 0 \\ b/a & 1\end{bmatrix} V_t$ converts \eqref{dV4} into
      \begin{equation} \label{dV5}
     d\widehat{V}_t  = \begin{bmatrix}-2\mu & 0 \\ 0 & 0 \end{bmatrix}\widehat{V}_t dt + \frac{\sigma}{\sqrt{\mu a}} \begin{bmatrix} -b & 1 \\ -a^2-b^2  & b \end{bmatrix} \widehat{V}_t dW_t^\phi
     \end{equation}
By Theorem 1 of Auslender and Milstein \cite{AM82}, the top Lyapunov exponent of \eqref{dV5} is $0 - \sigma^2 b^2/(2 \mu a) + O(\sigma^4)$ as $\sigma \to 0$.  Thus the top Lyapunov exponent of \eqref{dV3} is 
    $$
    -\frac{ \sigma^2 a}{2\mu} -\frac{\sigma^2 b^2}{2 \mu a} + O(\sigma^4) = - \frac{\sigma^2(a^2+b^2)}{2 \mu a}  + O(\sigma^4).
    $$ 
However this does not say anything directly about $\lambda(\mu,a,b,\sigma)$.  There is no proof here that approximation $r_t \sim \sqrt{\mu/a}$ is sufficiently accurate to ensure that \eqref{dV2} and \eqref{dV3} have the same top Lyapunov exponent up to $O(\sigma^4)$.  

For Theorem \ref{thm sigma0} the heuristic argument is less indicative of the actual method of proof.  Instead of adapting the method of \cite{AM82}, the proof of Theorem \ref{thm sigma0}, given in Section \ref{sec sigma0 proof}, uses the adjoint method to evaluate the Furstenberg-Khasminskii formula \eqref{lam khas} in the more typical setting of an asymptotic expansion for small $\sigma$. 

\section{Proof of Theorem \ref{thm binf}} \label{sec binf proof}   

\subsection{Tangent substitution}  \label{sec tangent}

Following the methods of Imkeller and Lederer \cite{IL99,IL01} we reparametrize the set $\R/(\pi \Z)$ using the new variable $z = \tan \psi$.  Thus $z$ takes values in the set $\widehat{\R}$ consisting of the extended real line $[-\infty,\infty]$ with the points $-\infty$ and $\infty$ identified.  The advantage of this reparametrization is that trigonometric functions in $\psi$ are replaced by rational functions of $z$.    Writing $z_t = \tan \psi_t$ we have
 \begin{align}
   dz_t  & =  \sec^2 \psi_t \, d\psi_t + \frac{1}{2} (2 \sec\psi_t)(\sec \psi_t \tan \psi_t)(d \psi_t)^2 \nonumber \\
       & =  \sec^2 \psi_t \left( 2r_t^2\cos \psi_t(b \cos \psi_t+a \sin \psi_t )dt-\frac{\sigma}{r_t}dW_t^\phi\right) + \sec^2 \psi_t \tan \psi_t \left(\frac{\sigma^2}{r_t^2}\right)dt  \nonumber \\
       & =  \left(2r_t^2  (b+az_t) + \frac{\sigma^2}{r_t^2}(1+z_t^2)z_t \right)dt - \frac{\sigma}{r}(1+z_t^2)\,dW_t^\phi. \label{z}
  \end{align}  
The behavior of $z_t$ near $\pm \infty$ is determined by the behavior of $1/z_t$ near 0, and the SDE for $1/z_t$ can be obtained from \eqref{z} by It\^o's formula.   %Note that the behavior of a function $f(z)$ near $z = \pm \infty$ can be determined from the behavior of $f(1/z)$ near $z = 0$.

  We will abuse notation by writing %$M = (0,\infty) \times \widehat{\R}$ in place of $(0,\infty) \times \R/\pi \Z$ and 
  $Q(r,z)$ in place of $Q(r,\psi)$ and by using $\widetilde{\nu}$ to denote the unique invariant measure for the process $\{(r_t,z_t): t \ge 0\}$ with values in $(0,\infty) \times \widehat{\R}$ given by (\ref{r},\ref{z}).  Thus   
  $$
     \lambda(\mu,a,b,\sigma) = \int_M Q(r, z) d\widetilde{\nu}(r,z)
      $$
where now $M = (0,\infty) \times \widehat{\R}$ and 
  $$
  Q(r,z) =\mu - ar^2+r^2\Bigl(b \sin 2\psi -a(1+\cos 2 \psi)\Bigr) = \mu - ar^2 + 2r^2\left(\frac{bz}{1+z^2}-\frac{a}{1+z^2}\right).  
   $$
We will also abuse notation by letting ${\cal L}$ denote the generator for the process $\{(r_t,z_t): t \ge 0\}$.  It may be written ${\cal L} = {\cal L}_r+{\cal L}_z$ where
   $$
  {\cal L}_r =  \left(\mu r - ar^3+\frac{\sigma^2}{2r}\right)\frac{\partial}{\partial r}+ \frac{\sigma^2}{2} 
 \frac{\partial^2}{\partial r^2}
 $$
and 
  $${\cal L}_z =   \left(r^2 (2b+2az) + \frac{\sigma^2}{r^2}(1+z^2)z \right)\frac{\partial}{\partial z} + \frac{\sigma^2}{2r^2}(1+z^2)^2 \frac{\partial^2}{\partial z^2}. 
  $$

\subsection{First application of the adjoint method} \label{sec first adj} 
 Here our choice of the function $F$ and the variable $w = z/A(r)$ is based on the heuristic argument in Section \ref{sec h binf}.  

\s

\n{\bf Definition.}  For a function $A = A(r)$ to be chosen later, define  
   \begin{equation} \label{F}
   F(r,z) = \frac{1}{2}\log\left(\frac{1+z^2}{A^2+1+ z^2}\right).
   \end{equation}
Note that we can treat $A$ as a constant for the purpose of calculating ${\cal L}_zF(r,z)$.

\begin{lemma} \label{lem QLFz}   
    \begin{align}
     Q(r,z) - {\cal L}_zF(r,z) & =   \frac{2br^2 z}{A^2+1+z^2} -\frac{A^2\sigma^2(1+z^2)(A^2+1-z^2)}{2r^2(A^2+1+z^2)^2}
     \nonumber \\
     &   \quad + \mu - ar^2 - 2ar^2\frac{A^2+1}{A^2+1+z^2} . \label{QLFz}
                   \end{align}
\end{lemma}  

\n{\bf Proof.}  This is a direct calculation.  We omit the details. \qed

\s

Writing $w=z/A$ in the right side of \eqref{QLFz} gives
  \begin{align}
     Q(r,z) - {\cal L}_zF(r,z) & =  \frac{2br^2}{A} \frac{w}{(1+A^{-2}+w^2)} -\frac{\sigma^2 A^2}{r^2} \frac{(A^{-2}+w^2)(1+A^{-2}-w^2)}{2(1+ A^{-2}+w^2)^2}
     \nonumber \\
     &   \quad + \mu - ar^2 - 2ar^2 \frac{1+A^{-2}}{(1+A^{-2}+w^2)}. \label{QLFw}
                   \end{align}
At this point we make the choice that 
   $2br^2/A = \sigma^2 A^2/r^2$, so that 
  $$
   A = A(r) = \left(\frac{2b r^4}{\sigma^2}\right)^{1/3}.
   $$
Then $2br^2/A = \sigma^2 A^2/r^2 = (2b \sigma r)^{2/3}$, and equation \eqref{QLFw} becomes 
    \begin{align}
      Q(r,z) - {\cal L}_zF(r,z) & =  (2b\sigma r)^{2/3}\left(\frac{w}{1+A^{-2}+w^2}  -\frac{(A^{-2}+w^2)(1+A^{-2}-w^2)}{2(1+ A^{-2}+w^2)^2}\right) \nonumber \\
      & \qquad + \mu - ar^2 - 2ar^2 \frac{1+A^{-2}}{(1+A^{-2}+w^2)}. \label{QLFw2}
        \end{align}
   Note that $A$ is now an explicit function of $r$, and that $A \to \infty$ as $b \to \infty$ for each $r >0$ but $A \to 0$ as $r \to 0$ for each $b$. The next two results quantify the effect of the $A^{-2}$ terms in the right side of \eqref{QLFw2} %deal with the behavior of the right side of \eqref{QLFw2} as $A \to \infty$
   and bound the contribution of ${\cal L}_rF(r,z)$ towards ${\cal L}F(r,z)$.  

\s

\n{\bf Definition.}  For $w \in \widehat{\R}$, define 
   $$
   \Phi(w) = \frac{w}{1+w^2}  -\frac{w^2(1-w^2)}{2(1+ w^2)^2} .
   $$
Note that 
  $$
  |\Phi(w)| \le \frac{|w|}{1+w^2}  +\frac{w^2|1-w^2|}{2(1+ w^2)^2} \le \frac{1}{2} + \frac{1}{2} = 1.
  $$

\begin{lemma}  \label{lem QLF est} Assume $b >0$.  With $A(r) = (2b r^4/\sigma^2)^{1/3}$ and $w =z/A(r)$ we have  
   \begin{equation} \label{QLF est}
   \left| Q(r,z) - {\cal L}_zF(r,z) - (2b\sigma r)^{2/3}\Phi(w) +\frac{2ar^2}{1+w^2} \right|   \le   \frac{3 \sigma (2b\sigma)^{1/3}}{4 r^{2/3}}+ \frac{a \sigma r^{2/3}}{(2b\sigma)^{1/3}} +  |\mu-ar^2|.
 \end{equation}
\end{lemma}

\n{\bf Proof.} Notice first that $1+\delta +w^2 \ge 1+\delta \ge 2 \sqrt{\delta}$, so that $\delta/(1+\delta+w^2) \le \sqrt{\delta}/2$.  Then 
   \begin{align*}
   \left| \frac{(\delta +w^2)(1+\delta-w^2)}{2(1+ \delta+w^2)^2} - 
  \frac{w^2(1-w^2)}{2(1+w^2)^2}\right|
     & = \frac{\delta(1+w^4+2w^6)+ \delta^2(1+w^2+2w^4)}{2(1+\delta+w^2)^2(1+w^2)^2}\\
     & \le \frac{\delta(1+w^2)^3+ \delta^2(1+w^2)^2}{(1+\delta+w^2)^2(1+w^2)^2}\\
     & =   \frac{\delta}{1+\delta+w^2} \le \frac{\sqrt{\delta}}{2}.
  \end{align*} 
Similarly
 $$
  \left| \frac{w}{1+\delta+w^2}   -   \frac{w}{1+w^2}\right|
   = \frac{\delta |w|}{(1+\delta+w^2)(1+w^2)} \le \frac{\delta}{2(1+\delta+w^2)} \le \frac{\sqrt{\delta}}{4}
   $$
and    
  $$
 \left|  \frac{1+\delta}{1+\delta+w^2}-\frac{1}{1+w^2}  \right| 
=  \left|  \frac{\delta w^2}{(1+\delta+w^2)(1+w^2)}  \right|  \le \frac{\delta}{1+\delta+w^2} \le \frac{\sqrt{\delta}}{2}.
$$  
Using these inequalities with $\delta = A(r)^{-2}% = \left(\dfrac{\sigma^2}{2br^4}\right)^{2/3}
$ in \eqref{QLFw2} gives 
  \begin{align*}
  \lefteqn{ \left| Q(r,z) - {\cal L}_zF(r,z) - (2b\sigma r)^{2/3}\Phi(w) +\frac{2ar^2}{1+w^2} \right|} \hspace{10ex} \\
  &   \le   (2b\sigma r)^{2/3}\left(\frac{1}{4A(r)} + \frac{1}{2A(r)}\right) + 2ar^2 \left(\frac{1}{2A(r)}\right)  +  |\mu-ar^2|,
  \end{align*}
and \eqref{QLF est} follows directly.  \qed

%\begin{remark} The estimate $\delta/(1+\delta + w^2) < \delta$ is stronger for small $\delta$, but verifying integrability for $r$ near 0, see Proposition \ref{prop lam}, involves small values of $A(r)$ and hence large values of $\delta$.  More precisely, with $\delta = A(r)^{-2}$ we see that $\sqrt{\delta}$ is $\nu$-integrable but $\delta$ is not $\nu$-integrable.
%\end{remark}

\begin{lemma} \label{lem LRF} With $A(r) = (2b r^4/\sigma^2)^{1/3}$ we have 
    \begin{equation} \label{LRF}
  |{\cal L}_r F(r,z)| \le \frac{4}{3}|\mu - ar^2| + \frac{8 \sigma (2b\sigma)^{1/3}}{9r^{2/3}}.
  \end{equation}
\end{lemma}

\n{\bf Proof.}  We have 
  \begin{align}
   {\cal L}_rF(r,z) & =  \left(\mu r - ar^3+\frac{\sigma^2}{2r}\right)\frac{\partial F}{\partial r}(r,z)+ \frac{\sigma^2}{2} 
\cdot  \frac{\partial^2 F}{\partial r^2}(r,z) \nonumber  \\[1ex]
    & =   -(\mu r - ar^3)\frac{A}{(A^2+1+z^2)} \frac{\partial A}{\partial r}-\frac{\sigma^2}{2r}\cdot\frac{A}{(A^2+1+z^2)} \frac{\partial A}{\partial r}  \nonumber \\
    &  \mbox{} \quad  + \frac{\sigma^2}{2}\cdot \frac{(A^2-1-z^2)}{(A^2+1+z^2)^2} \left(\frac{\partial A}{\partial r}\right)^2- \frac{\sigma^2}{2}\cdot \frac{A}{(A^2+1+z^2)} \frac{\partial^2 A}{\partial r^2}  \nonumber \\[1ex]
    & =  I + II + III + IV, \label{LRF sum}
     \end{align}
say.  Write $A = A(r) =k r^{4/3}$, where $k = (2b/\sigma^2)^{1/3}$.  Then $\dfrac{\partial A}{\partial r} = \dfrac{4}{3}kr^{1/3}$ and $\dfrac{\partial^2 A}{\partial r^2} = \dfrac{4}{9}kr^{-2/3}$.  In the following we use $\dfrac{\partial A}{\partial r} = \dfrac{4A}{3r}$ and $A^2+1+z^2 \ge A^2+1 \ge 2A$.  We have 
   $$
   |I| = |\mu r-ar^3|\cdot \frac{4A^2}{3r(A^2+1+z^2)} \le \frac{4}{3}|\mu - ar^2|
   $$ 
and    $$
   |II| = \frac{\sigma^2}{2r}\cdot \frac{4A^2}{3r(A^2+1 +z^2)} \le \frac{\sigma^2}{2r}\cdot \frac{2A}{3r}
   = \frac{\sigma(2b\sigma)^{1/3}}{3 r^{2/3}}
   $$ 
and
$$
   |III| = \frac{\sigma^2}{2}\cdot \frac{|A^2-1-z^2|}{(A^2+1+z^2)^2}  \cdot \frac{16A^2}{9r^2} \le \frac{\sigma^2}{2}\cdot \frac{8A}{9r^2 } =  \frac{4\sigma(2b\sigma)^{1/3}}{9 r^{2/3}}
     $$ 
and 
  $$
   |IV| = \frac{\sigma^2}{2} \cdot \frac{A}{(A^2+1+z^2)} \cdot \frac{4}{9}kr^{-2/3} \le \frac{\sigma^2}{2}\cdot \frac{1}{2} \cdot \frac{4}{9}kr^{-2/3}=  \frac{\sigma(2b\sigma)^{1/3}}{9 r^{2/3}}.
     $$ 
The result now follows from \eqref{LRF sum}.  
   \qed

\begin{proposition} \label{prop QLF}   Assume $b >0$.  With $A(r) = (2b r^4/\sigma^2)^{1/3}$
 \begin{align}  \nonumber   
   \lefteqn{  \left| Q(r,z) - {\cal L}F(r,z) - (2b\sigma r)^{2/3}\Phi\left(\frac{z}{A(r)}\right)+ \frac{2ar^2}{1+(z/A(r))^2}\right|} \hspace{30ex}\\[2ex] 
   & \le \frac{59 \sigma(2b \sigma)^{1/3}}{36r^{2/3}} + \frac{a \sigma r^{2/3}}{(2b\sigma)^{1/3}} + \frac{7}{3}|\mu-ar^2|.\label{QLF}
    \end{align}
\end{proposition}      

\n{\bf Proof.}  This is an immediate consequence of Lemmas \ref{lem QLF est} and \ref{lem LRF}.  \qed

\begin{corollary} \label{cor QLF}   Assume $b >0$.  With $A(r) = (2b r^4/\sigma^2)^{1/3}$ we have   
   \begin{align} \nonumber  
    \lefteqn{\left|\lambda(\mu,a,b,\sigma) - (2b\sigma)^{2/3}\int_M r ^{2/3} \Phi\left(\frac{z}{A(r)}\right)d\widetilde{\nu}(r,z)+ 2a\int_M \frac{r^2}{1+(z/A(r))^2}d\widetilde{\nu}(r,z)\right|}  \\
    & \le  \frac{59 \sigma (2b\sigma)^{1/3}}{36} \int_{\R_+} r^{-2/3}d\nu(r) + \frac{a  \sigma}{(2b\sigma)^{1/3}}\int_{\R_+} r^{2/3}d\nu(r) +\frac{7}{3}\int_{\R_+}|\mu-ar^2|d\nu(r). \label{lam est}
    \end{align}        
\end{corollary}

\n{\bf Proof.}  
From the formula \eqref{rho} for the density of $\nu$ we have $\int_{\R_+} r^\kappa d\nu(r) < \infty$ for all $\kappa > -2$.  Therefore all of the terms in \eqref{QLF} are integrable with respect to $\widetilde{\nu}$, and so 
   \begin{align*}
  \lefteqn{\left|
  \lambda(\mu,a,b,\sigma) - (2b\sigma )^{2/3}\int_M r^{2/3}\Phi\left(\frac{z}{A(r)}\right)d\widetilde{\nu}(r,z) +\int_M \frac{2ar^2}{1+(z/A(r))^2}d\widetilde{\nu}(r,z)\right|} \hspace{4ex}\\
 & =    \left| \int_M \left(Q(r,z) - (2b\sigma r)^{2/3}\Phi\left(\frac{z}{A(r)}\right)+ \frac{2ar^2}{1+(z/A(r))^2}\right)d\widetilde{\nu}(r,z)\right| \\
   & \le     \left|\int_M {\cal L}F(r,z)d\widetilde{\nu}(r,z)\right| + \frac{59\sigma (2 b \sigma)^{1/3}}{36} \int_{\R_+} r^{-2/3} d\nu(r) + \frac{a\sigma}{(2b\sigma)^{1/3}}\int_{\R_+} r^{2/3}d\nu(r)\\
   & \qquad  +\frac{7}{3}\int_{\R_+} |\mu-ar^2|d\nu(r).
      \end{align*}
It remains to show that $\int_M {\cal L}F(r,z)d\widetilde{\nu}(r,z) = 0$.    Notice first that $-\frac{1}{2}\log (1+A(r)^2) \le F(r,z) \le 0$ so that $|F(r,z)| \le \frac{1}{2}\log \bigl(1+(2b/\sigma^2)^{2/3} r^{8/3}\bigr)$.  Choose positive $V \in C^2(\R_+)$ so that $V(r)  =1$ for $0 < r \le 1$ and $V(r) =r$ for $r \ge 2$.  Then $V$ is $\widetilde{\nu}$ integrable and $\sup\{{\cal L}V(r)/V(r): r > 0\} < \infty$ and $|F(r,z)|/V(r) \to 0$ as $r \to \infty$.  %Moreover the estimate \eqref{QLF} shows there exists $\alpha > 1$ such that $|{\cal L}F(r,z)|^\alpha$ is $\widetilde{\nu}$ integrable. 
 Therefore we may apply Proposition \ref{prop-bg} from the Appendix and deduce that  $\int_M {\cal L}F(r,z)d\widetilde{\nu}(r,z) = 0$.  \qed 

\s

Now consider the variable $w = z/A(r)$ where $A(r) = (2br^4/\sigma^2)^{1/3}$, and the corresponding process $w_t = z_t/A(r_t)$.  The mapping $(r,z) \mapsto (r,z/A(r)) = (r,w)$ is a diffeomorphism of $M = (0,\infty) \times \widehat{\R}$.  It converts the diffusion $\{(r_t,z_t): t \ge 0\}$ into a diffusion $\{(r_t,w_t): t \ge 0\}$ with generator ${\cal M}$, say.  
 Let $\eta$ denote the pushforward of $\widetilde{\nu}$ under the diffeomorphism.  Then $\eta$ is the unique invariant probability measure for the diffusion $\{(r_t,w_t): t \ge 0\}$. % given by (\ref{r},\ref{w}).  
 Note that the $r$ marginal measure of $\eta$ is $\nu$. % as $\widetilde{\nu}$.  %Let ${\cal M}$ denote the generator for the system $\{(r_t,w_t): t \ge 0\}$.  

\begin{lemma} \label{lem w}  Assume $b >0$.  With $A(r) = (2b r^4/\sigma^2)^{1/3}$ and $w_t = z_t/A(r_t)$ we have
    \begin{align}
  dw_t  & =  (2b\sigma r_t)^{2/3}(1+w_t^3)dt- (2b\sigma r_t)^{1/3} w_t^2 dW_t^\phi \nonumber \\
      &  \qquad + \left( -\frac{4\mu}{3}+\frac{10ar_t^2}{3} + \frac{17 \sigma^2}{9r_t^2}\right)w_t dt - \frac{4\sigma w_t}{3r_t}dW_t^r - \frac{\sigma^2}{(2b\sigma r_t)^{1/3}r_t^2}dW_t^\phi. \label{w} 
  \end{align} 
\end{lemma}

\n{\bf Proof.}  Write  $A_t = A(r_t) = kr_t^{4/3}$ where $k = (2b/\sigma^2)^{1/3}$.   We have  
  \begin{align*}
   dA_t  =  \frac{4k}{3} r_t^{1/3}dr_t+ \frac{1}{2}\cdot \frac{4k}{9} r_t^{-2/3}(dr_t)^2
        & = \frac{4A_t}{3r_t} dr_t+ \frac{2A_t}{9 r_t^2}(dr_t)^2 \\
       & = A_t \left(  \frac{4}{3}(\mu  - ar_t^2)+\frac{8 \sigma^2}{9r_t^2}\right)dt+ \frac{4 \sigma A_t}{3r_t} dW_t^r.
        \end{align*}
%OLD We have  
 % \begin{align*}
 %  \frac{1}{A_t}dA_t & =  \frac{1}{A_t}\left(\frac{4k}{3} r_t^{1/3}dr_t+ \frac{1}{2}\cdot \frac{4k}{9} r_t^{-2/3}(dr_t)^2\right) \\
 %      & = \left(  \frac{4}{3}(\mu  - ar_t^2)+\frac{8 \sigma^2}{9r_t^2}\right)dt+ \frac{4 \sigma}{3r_t} dW_t^r.
 %       \end{align*}
Recall equation \eqref{z} for $\{z_t: t \ge 0\}$.  Since $\{W_t^r: t \ge 0\}$ and $\{W_t^\phi: t \ge 0\}$ are independent we have $dz_tdA_t = 0$. Therefore  
  \begin{align*}
  dw_t  & =  \frac{1}{A_t}dz_t - \frac{z_t}{A_t^2}dA_t + \frac{z_t}{A_t^3}(dA_t)^2 \\
  & =  \frac{1}{A_t}\left(2r_t^2 (b+az_t) + \frac{\sigma^2}{r_t^2}(1+z_t^2)z_t \right)dt - \frac{\sigma}{A_t r_t}(1+z_t^2)\,dW_t^\phi   - \frac{z_t}{A_t^2}dA_t + \frac{z_t}{A_t^3}(dA_t)^2\\
     & =  \frac{1}{A_t}\left(2r_t^2 (b+aA_tw_t) + \frac{\sigma^2}{r_t^2}(1+A_t^2w_t^2)A_t w_t \right)dt - \frac{\sigma}{A_t r_t}(1+A_t^2 w_t^2)\,dW_t^\phi \\
     & \qquad  - \frac{w_t}{A_t}dA_t + \frac{w_t}{A_t^2}(dA_t)^2\\
      & =  \left( \frac{2b r_t^2}{A_t}+2ar_t^2w_t + \frac{\sigma^2 w_t}{r_t^2} + \frac{\sigma^2 A_t^2 w_t^3}{r_t^2}   \right)dt - \frac{\sigma }{r_t}\left(\frac{1}{A_t}+A_tw_t^2\right)dW_t^\phi\\
       & \qquad  - w_t \left(  \frac{4}{3}(\mu  - ar_t^2)+\frac{8\sigma^2}{9r_t^2}\right)dt - \frac{4 \sigma w_t}{3r_t} dW_t^r + \frac{16 \sigma^2 w_t}{9r_t^2}dt, 
  \end{align*}  
and \eqref{w} follows by putting $A_t = (2b \sigma r_t^4/\sigma^2)^{1/3}$.   \qed

\begin{corollary} \label{cor lam} Assume $b >0$.  We have  
   \begin{align} \nonumber  
    \lefteqn{\left|\lambda(\mu,a,b,\sigma) - (2b\sigma)^{2/3}\int_M r ^{2/3} \Phi(w)d\eta(r,w)+ 2a\int_M \frac{r^2}{1+w^2}d\eta(r,w)\right|} \\
    & \le \frac{59 \sigma (2b\sigma)^{1/3}}{36} \int_{\R_+} r^{-2/3}d\nu(r) + \frac{a  \sigma}{(2b\sigma)^{1/3}}\int_{\R_+} r^{2/3}d\nu(r) +\frac{7}{3}\int_{\R_+}|\mu-ar^2|d\nu(r), \label{lam est2}
    \end{align}        
where $\eta$ is the unique invariant probability measure for the diffusion $\{(r_t,w_t): t \ge 0\}$ on $M = (0,\infty) \times \widehat{R}$ with generator 
   \begin{equation} \label{M}
   {\cal M} = (2b\sigma r)^{2/3}{\cal M}_0 + {\cal M}_1+ (2b \sigma r)^{-2/3} {\cal M}_2
   \end{equation}
where 
    \begin{align*}
  {\cal M}_0 & = (1+ w^3)\frac{\partial}{\partial w} + \frac{w^4}{2}\frac{\partial^2}{\partial w^2}, \\
 {\cal M}_1 & =  \left(\mu r_t - ar_t^3+\frac{\sigma^2}{2r_t}\right)\frac{\partial}{\partial r} + 
       \left( -\frac{4\mu}{3}+\frac{10ar^2}{3} + \frac{17 \sigma^2}{9r^2}\right)w \frac{\partial}{\partial w}\\
       &  \mbox{} \qquad
    +\frac{\sigma^2}{2}\frac{\partial^2}{\partial r^2} - \frac{4\sigma^2 w}{3r}\frac{\partial ^2}{\partial r\partial w} + \frac{17 \sigma^2w^2}{9r^2} \frac{\partial ^2}{\partial w^2}, \\
    {\cal M}_2 & = \frac{\sigma^4}{2r^4}\frac{\partial^2}{\partial w^2}.
  \end{align*}
  
  \end{corollary}

\n{\bf Proof.}  The inequality \eqref{lam est2} follows immediately from \eqref{lam est} because $\eta$ is the pushforward of $\widetilde{\nu}$.  The formula for ${\cal M}$ follows from the equations (\ref{r},\ref{w}) for $r_t$ and $w_t$.  \qed

\subsection{$b \to \infty$ with fixed $\sigma > 0$}  Corollary \ref{cor lam} is designed to be useful for the proofs of both Theorem \ref{thm binf} and Theorem \ref{thm CE unif}.  In this section we assume $\sigma > 0$ is fixed as $b \to \infty$, and so we concentrate on the term $(2 b \sigma)^{2/3}\int _M r^{2/3}\Phi(w) d\eta(r,w)$ in the left side of \eqref{lam est2}.  The other term $2a \int_M \frac{r^2}{1+w^2}d\eta(r,w)$ will merely contribute to the remainder in this section, although it plays an important role in the proof of Theorem \ref{thm CE unif} in Section \ref{sec CE proof}.  From the expression \eqref{M} for the generator ${\cal M}$ we see that the highest order term ${\cal M}_0$ does not involve the parameter $r$.  This corresponds to the heuristic argument in Section \ref{sec h binf} about separation of time scales.  Technically it means that it will be useful to consider an adjoint problem involving the function $\Phi(w)$ and the generator ${\cal M}_0$ on the space $\widehat{\R}$

\begin{lemma} \label{lem G}  There is a smooth bounded function $G: \widehat{\R} \to \R$ and a constant $\gamma_0 = \pi/ \bigl(2^{1/3} 3^{1/6}[\Gamma(1/3)]^2 \bigr)\approx 0.29$ such that 
     \begin{equation} \label{G}
     {\cal M}_0 G(w)
     = \Phi(w) - \gamma_0.
   \end{equation}
Moreover $(1+w^2) G'(w)$ and $(1+w^2) G''(w)$ are bounded.
 \end{lemma} 
 
 \n{\bf Proof.}   Consider the linear SDE \eqref{YA}:
  $$
   dY_t = \begin{bmatrix} 0 & 0 \\1 & 0 \end{bmatrix}Y_tdt
   +\begin{bmatrix} 0 & 1 \\0 & 0 \end{bmatrix}Y_tdW_t.
   $$
Writing $Y_t = \|Y_t\| \begin{bmatrix} \cos \vartheta_t \\ \sin \vartheta_t \end{bmatrix}$ and then $w_t = \tan \vartheta_t$ we have 
    \begin{equation} \label{dwA}
    dw_t = (1+w_t^3)dt - w_t^2 dW_t
    \end{equation}
and 
   $$
   d\log \|Y_t\| = \left(\frac{w_t}{1+w_t^2} - \frac{w_t^2(1-w_t^2)}{2(1+w_t^2)^2}\right)dt + \frac{w_t}{1+w_t^2}dW_t = \Phi(w_t)dt +  \frac{w_t}{1+w_t^2}dW_t .
   $$  
Thus ${\cal M}_0$ is the generator for the diffusion $\{w_t: t \ge 0\}$ given by \eqref{dwA}, and $\Phi(w)$ is the integrand in the Furstenberg-Khasminskii formula for the top Lyapunov exponent $\gamma_0$ of \eqref{YA}.  Therefore ${\cal M}_0 G(w) = \Phi(w) - \gamma_0
    $ is the adjoint equation for \eqref{YA}. %, and $\gamma_0$ is its top Lyapunov exponent.  
 Theorem 5.1 of Arnold, Oeljeklaus and Pardoux \cite{AOP86} applies to the linear SDE \eqref{YA} and gives the existence of the function $G$, and Ariaratnam and Xie \cite{AX90} gives the exact value of $\gamma_0$.  Moreover there is a smooth $\pi$-periodic function $\widetilde{G}$ such that $\widetilde{G}(\vartheta) = G(\tan \vartheta)$.   Therefore $G$ is bounded.  Also
   $$
   \widetilde{G}'(\vartheta)  = (\sec^2 \vartheta )G'(\tan \vartheta) = (1+w^2)G'(w).
   $$
Since $\widetilde{G}'(\vartheta)$ is bounded, then so is $(1+w^2)G'(w)$.  Moreover  
    $$
   \widetilde{G}''(\vartheta)  =( 2 \sec^2 \vartheta \tan \vartheta) G'(\tan \vartheta) + (\sec^4 \vartheta) G''(\tan \vartheta) = (1+w^2)\Bigl(2w G'(w)+(1+w^2)G''(w)\Bigr)
   $$  
so that $2w G'(w)+(1+w^2)G''(w)$ is bounded.  Since $|2wG'(w)| \le (1+w^2|G'(w)|$ is bounded it follows that $(1+w^2)G''(w)$ is bounded.   \qed
   
\s

%Now 
%  \begin{eqnarray*}
%   {\cal M} G(r,w) & = & (2b\sigma r)^{2/3}{\cal M}_0G(w) + {\cal M}_1G(r,w)+
%      (2b \sigma r)^{-2/3} {\cal M}_2G(r,w)\\
%      & = & (2b\sigma r)^{2/3}(\Phi(w)-\lambda_0) + {\cal M}_1G(r,w)+
%      (2b \sigma r)^{-2/3} {\cal M}_2G(r,w)
% \end{eqnarray*}
%where the terms ${\cal M}_1G$ and $(2b \sigma r)^{-2/3}{\cal M}_2G$ involve products of the form $f(r)g(w)$ where each $g(w)$ is bounded but where some of the functions $f(r)$ are not $\nu$ integrable.    The following proof deals with this lack of integrability.    

\begin{proposition}  \label{prop Phi}   Given $\mu \in \R$, $a > 0$ and $\sigma >0$,   
 $$
 \int_M r^{2/3}\Phi(w)d\eta(r,w) \to  \gamma_0 \int_{\R_+} r^{2/3}d\nu(r)
 $$
as $b \to \infty$.
  \end{proposition}

%\begin{remark}  Need to restate this in a more uniform way.  Possibly
%   $$
%   \left| (2b\sigma)^{2/3}\int r^{2/3}\Phi(w)d\eta(r,w) - (2b\sigma)^{2/3}\gamma_0 \int r^{2/3}d\widetilde{\nu}(r) %\right| \to 0
%   $$
%as $b \sigma \to \infty$, uniformly for ???
%   \end{remark}     
\n{\bf Proof.}  Using Corollary \ref{cor lam} and Lemma \ref{lem G} we have  
  \begin{align*}
   {\cal M} G(r,w) & =  (2b\sigma r)^{2/3}{\cal M}_0G(w) + {\cal M}_1G(r,w)+
      (2b \sigma r)^{-2/3} {\cal M}_2G(r,w)\\
      & =  (2b\sigma r)^{2/3}(\Phi(w)-\lambda_0) + {\cal M}_1G(r,w)+
      (2b \sigma r)^{-2/3} {\cal M}_2G(r,w)
 \end{align*}
where the terms ${\cal M}_1G$ and $r^{-2/3}{\cal M}_2G$ involve products of the form $f(r)g(w)$ where each $g(w)$ is bounded but where some of the functions $f(r)$ are not $\nu$ integrable.  The following proof uses truncation to avoid this lack of integrability. Note that the measure $\eta$ depends on $b$, but its marginal $\nu$ does not depend on $b$.  

Let $C \subset (0,\infty)$ be a closed compact interval and let $h:(0,\infty) \to [0,1]$ be a smooth compactly supported function with $h\big|_{C} = 1$.  Define $G_h(r,w) = h(r)G(w)$. Since ${\cal M}_0$ does not involve derivatives with respect to $r$ we have ${\cal M}_0G_h(r,w) = h(r){\cal M}_0G(w) = h(r)(\Phi(w)-\gamma_0)$.   
     Then
 \begin{equation} \label{MGh}
   {\cal M} G_h(r,w)  =  (2b\sigma r)^{2/3}h(r)(\Phi(w)-\lambda_0) + {\cal M}_1G_h(r,w)+
      (2b \sigma r)^{-2/3} {\cal M}_2G_h(r,w).
 \end{equation} 
The terms ${\cal M}_1G_h(r,w)$ and $r^{-2/3}{\cal M}_2G_h(r,w)$ involve a finite sum of products of the form $cf(r)g(w)$ where each $g(w)$ is $G(w)$ or $wG'(w)$ or $G''(w)$ or $w^2G(w)$ and hence bounded, and each $f(r)$ is a power of $r$ times a multiple of $h(r)$ or $h'(r)$ or $h''(r)$ and hence bounded, and each constant $c$ is a multiple of one of the parameters $\mu$ or $a$ or $\sigma^2$ or $\sigma^4$.  Therefore ${\cal M}_1G_h(r,w)$ and $r^{-2/3}{\cal M}_2G_h(r,w)$ are both bounded.  Moreover since $G$ and $\Phi$ are bounded, it follows (using \eqref{MGh}) that $G_h$ and ${\cal M}G_h$ are bounded.        
 Using Proposition \ref{prop-bg} from the Appendix we have 
   \begin{align*}
     0 & =  \int_M {\cal M}G_h(r,w)d\eta(r,w) \\
    & =  \int_M (2b\sigma r)^{2/3} h(r) \left(\Phi(w) - \gamma_0\right) d\eta(r,w) + \int_M {\cal M}_1G_h(r,w)d\eta(r,w)\\
    &  \qquad +
     \int_M  (2b \sigma r)^{-2/3} {\cal M}_2G_h(r,w)d\eta(r,w).
     \end{align*} 
Now $|\Phi(w) - \gamma_0|$ is bounded, by $c$, say.  Then
   \begin{align}
  \lefteqn{ \left|\int_M r^{2/3}\Phi(w)d\eta(r,w) - \gamma_0 \int_{\R_+} r^{2/3}d\nu(r) \right| } \hspace{5ex} \nonumber \\
  % & =  \left|\int_M r^{2/3}(\Phi(w)-\gamma_0)d\eta(r,w)\right| \nonumber\\
   & \le   \left|\int_M r^{2/3}(1-h(r))(\Phi(w)-\gamma_0)d\eta(r,w)\right| +  \left|\int_M r^{2/3}h(r)(\Phi(w)-\gamma_0)d\eta(r,w)\right|   \nonumber \\
   & \le   c \int_{\R_+ \setminus C} r^{2/3}d\nu(r) + (2b\sigma)^{-2/3}\left|\int_M {\cal M}_1G_h(r,w)d\eta(r,w)\right| \nonumber \\
   &  \qquad +  (2b\sigma)^{-4/3}\left|\int_M r^{-2/3} {\cal M}_2G_h(r,w)d\eta(r,w)\right| \nonumber \\
    & \le   c \int_{\R_+ \setminus C} r^{2/3}d\nu(r) + (2b\sigma)^{-2/3} \sup_{(r,w)}\Bigl|{\cal M}_1G_h(r,w)\Bigr| \nonumber \\
   &  \qquad 
   +  (2b\sigma)^{-4/3}\sup_{(r,w)}\Bigl| r^{-2/3} {\cal M}_2G_h(r,w)\Bigr|. \label{Phi est}
   \end{align}
Therefore
   $$
   \limsup_{b \to\infty}   \left|\int r^{2/3}\Phi(w)d\eta(r,w) - \gamma_0 \int r^{2/3}d\nu(r) \right|\le c \int_{\R_+ \setminus C} r^{2/3}d\nu(r) 
   $$
and the proof is completed by letting $C \nearrow \R_+$.  \qed

%\begin{theorem} \label{thm binfA}  Given $\mu \in\R$, $a >0$ and $\sigma >0$,     
% \begin{equation} \label{binfA}
%       \lim_{b \to \infty} \frac{\lambda(\mu,a,b,\sigma)}{(2b\sigma)^{2/3}} = \gamma_0 \int_{\R_+} r^{2/3}d\nu(r).
%    \end{equation}
%\end{theorem}    

\s

\n{\bf Proof of Theorem \ref{thm binf}.}    We have 
   \begin{align}
  \lefteqn{ \left|\frac{\lambda(\mu,a,b,\sigma)}{(2b\sigma)^{2/3}} - \gamma_0 \int_{\R_+} r^{2/3}d\nu(r) \right|} \hspace{8ex} \nonumber \\
      & \le  
   \left|\frac{\lambda(\mu,a,b,\sigma)}{(2b\sigma)^{2/3}} - \int_M r^{2/3}\Phi(w)d\eta(r,w)
     + \frac{2a}{(2b\sigma)^{2/3}}\int_M \frac{r^2}{1+w^2}d\eta(r,w) \right|  \nonumber \\
      &  \qquad + \frac{2a}{(2b\sigma)^{2/3}}\int_{\R_+} r^2 d\nu(r)+ 
      \left| \int_M r^{2/3}\left(\Phi(w)-\gamma_0\right)d\eta(r,w) \right|,   \label{triangle}
\end{align}
and the result follows by using Corollary \ref{cor lam} and Proposition \ref{prop Phi} on the terms on the right.  \qed 

\section{Proofs for the uniform versions of Theorem \ref{thm binf}} \label{sec proof unif}

Recall we write $\nu = \nu_{(\mu,a,\sigma)}$ to denote the dependence of the $r_t$-invariant probability measure on the parameters $\mu \in \R$, $a> 0$ and $\sigma> 0$.

\s

\n{\bf Proof of Theorem \ref{thm binf unif}.}  (i) This follows from a detailed analysis of the inequality \eqref{triangle} and in particular the terms in the inequalities \eqref{lam est2} and \eqref{Phi est}.  For the inequality \eqref{lam est2} the Cauchy-Schwarz inequality 
    $$
     \int_{\R_+} r^{2/3}d\nu_{(\mu,a,\sigma)}(r) \le \left(\int_{\R_+} r^{-2/3}d\nu_{(\mu,a,\sigma)}(r)\right)^{1/2}\left(\int_{\R_+} r^2 d\nu_{(\mu,a,\sigma)}(r)\right)^{1/2}
    $$
implies that $\int_{\R_+} r^{2/3}d\nu_{(\mu,a,\sigma)}(r)$ is bounded above for $(\mu,a,\sigma) \in P$.  For the inequality \eqref{Phi est}, if $C =[\kappa_1,\kappa_2] \subset (0,\infty)$ then
    \begin{align*}
    \int_{\R_+ \setminus C} r^{2/3} d\nu_{(\mu,a,\sigma)}(r) 
   &  = \int_{(0,\kappa_1)} r^{2/3} d\nu_{(\mu,a,\sigma)}(r)  + \int_{(\kappa_2,\infty)} r^{2/3} d\nu_{(\mu,a,\sigma)}(r)\\
   & \le k_1^{2/3} + \kappa_2^{-4/3} \int_{\R_+} r^2 d\nu_{(\mu,a,\sigma)}(r)  
    \end{align*}  
so that $ \int_{\R_+ \setminus C} r^{2/3} d\nu_{(\mu,a,\sigma)}(r) \to 0$ as $C \nearrow \R_+$ uniformly for $(\mu,a,\sigma) \in P$.   Finally, since $P$ is bounded then given $C$ and hence $h$ the functions ${\cal M}_1G_h$ and ${\cal M}_2G_h$ are bounded uniformly for $(\mu,a,\sigma) \in P$.  

(ii) This follows from (i) by choosing $\e < \gamma_0 \inf\{\int_{\R_+} r^{2/3} d\nu_{(\mu,a,\sigma)}(r): (\mu,a,\sigma) \in P\}$.  \qed

\s

%For applications of Theorem \ref{thm binfA unif} we note that the comparison theorem, , for the process $\{r_t: t \ge 0\}$ with different values for $\mu$ implies a corresponding ordering of the invariant probabilities $\nu_{(\mu,a,\sigma)}$.  Thus if $\mu_1 < \mu_t$ then $\int_{\R_+} f(r)d\nu_{(\mu_1,a,\sigma)}(r) \le \int_{\R_+} f(r)d\nu_{(\mu_2,a,\sigma)}(r)$ for all increasing functions $f:\R_+ \to \R$. 

%Similarly for $a$ with changes in sign.
   
%\begin{corollary} \label{cor binfA unif}  (i)  Fix $a > 0$ and $\sigma > 0$.  Given $\mu_0 > 0$ there exists $k_1$ such that $\lambda(\mu,a,b,\sigma) >0$ if $-\mu_0 \le \mu \le \mu_0$ and $b \ge k_1$.

%(ii) Fix $\mu > 0$ and $a > 0$.  Given $\sigma_0 > 0$ there exists $k_2$ such that $\lambda(\mu,a,b,\sigma) >0$ if $0 < \sigma \le \sigma_0$ and $b\sigma  \ge k_2$.

%(iii) %Fix $\mu < 0$ and $\sigma > 0$.  Given $a_0 > 0$ there exists $k_3$ such that $\lambda(\mu,a,b,\sigma) >0$ if $0 < a \le a_0$ and $b \ge k_3$.
%   Fix $\mu <0$ and $a > 0$.  Given $\sigma_0 > 0$ there exists $k_3$ such that $\lambda(\mu,a,b,\sigma) > 0$ if $0 <\sigma \le \sigma_0$ and $b\sigma^2 \ge k_3$.
%   \end{corollary}
   
\n{\bf Proof of Corollary \ref{cor binf unif}.}   (i)  Fix $a > 0$ and $\sigma > 0$ and $\mu_0 >0$.  The comparison theorem for one-dimensional diffusions, see \cite[Theorem 1]{IK77}, implies the diffusion $\{r_t: t \ge 0\}$ given by \eqref{r} is stochastically monotone in the parameter $\mu$.  It follows that the invariant measures $\nu_{(\mu,a,\sigma)}$ are stochastically monotone in $\mu$.  That is, if $\mu_1 \le \mu_2$ then $\int_{\R_+} f(r)d\nu_{(\mu_1,a,\sigma)}(r) \le \int_{\R_+} f(r)d\nu_{(\mu_2,a,\sigma)}(r)$ for all increasing functions $f:\R_+ \to \R$.   The stochastic monotonicity implies for $-\mu_0 \le \mu \le \mu_0$ we have 
    $$
   0 <  \int r^\gamma d\nu_{-\mu_0,a,\sigma}(r) \le  \int r^\gamma d\nu_{\mu,a,\sigma}(r) \le  \int r^\gamma d\nu_{\mu_0,a,\sigma}(r) < \infty, \qquad \gamma > 0
    $$
and   
   $$
   0 <  \int r^\gamma d\nu_{\mu_0,a,\sigma}(r) \le  \int r^\gamma d\nu_{\mu,a,\sigma}(r) \le  \int r^\gamma d\nu_{-\mu_0,a,\sigma}(r) < \infty, \qquad  -2 < \gamma < 0.    $$
Therefore the set $P =\{(\mu,a,\sigma): -\mu_0 \le \mu \le \mu_0\}$ satisfies the assumptions for both parts of Theorem \ref{thm binf unif}.  Since $\sigma >0$ is fixed, the condition $b\sigma \ge k$ is equivalent to the condition $b \ge k/\sigma = k_1$, say.

(ii)  Fix $\mu > 0$ and $a > 0$ and $\sigma_0 > 0$.  Lemma \ref{lem nu bounds} in the Appendix shows that the set $P =\{(\mu,a,\sigma): 0 < \sigma \le \sigma_0\}$ satisfies the assumptions for both parts of Theorem \ref{thm binf unif}.  

(iii)  Note first that when $\mu < 0$ the SDE \eqref{r} for the diffusion $\{r_t: t \ge 0\}$ is well defined and positive recurrent when $a=0$, so that $\nu_{(\mu,0,\sigma)}$ is a well-defined probability measure.  Now fix $\mu < 0$ and $a > 0$ and and $\sigma_0 > 0$.  The measures $\nu_{(\mu,a\sigma^2,1)}$ are stochastically decreasing for $0 \le \sigma \le \sigma_0$, so the method used in (i) implies that the set $P = \{(\mu,a \sigma^2,1): 0 < \sigma \le \sigma_0^2\}$ satisfies the assumptions for both parts of Theorem \ref{thm binf unif}. So there exists $k_3$ such that $\lambda(\mu,a\sigma^2,\tilde{b},1) > 0$ if $0 < \sigma  \le \sigma_0^2$ and $\tilde{b} \ge k_3$.   Finally use the scaling result
      $$
      \lambda(\mu,a,b,\sigma) = \lambda(\mu,a\sigma^2,b\sigma^2,1)
      $$
with $\tilde{b} = b\sigma^2$.  \qed

\section{Proof of Theorem \ref{thm CE unif}.} \label{sec CE proof}
   
In this section $b\sigma$ is assumed to be in a compact set $D \subset (0,\infty)$, so that it is bounded away from 0 and $\infty$, while $\sigma$ is small.   We return to Corollary \ref{cor lam} and see that the two integrals in the left side of \eqref{lam est2} are of comparable size.  Thus we concentrate on the term  
  \begin{equation} \label{Phi r int}
   \int_M \left((2b\sigma)^{2/3} r ^{2/3} \Phi(w) - \frac{2a r^2}{1+w^2}\right)d\eta(r,w)
  \end{equation}
 where $\eta$ is the invariant probability measure for the process $\{(r_t,w_t): t \ge 0\}$ with generator ${\cal M}$.  
 For $\mu > 0$ and small $\sigma$ the distribution of $r$ is concentrated near $\sqrt{\mu/a} \equiv r_0$.  Lemma \ref{lem Phi zeta} deals with the effect of replacing $r$ with $r_0$ in the integrand of \eqref{Phi r int}, and Lemma \ref{lem M CE} shows the new hierarchy of larger and smaller terms in the formula \eqref{M} for ${\cal M}$ when $\sigma$ is small and $r$ is close in distribution to $r_0$.  % in oth The next two lemmas show the effect of this on the integrand in  and on the hierarchy  in the feffect ormula . 

\s

\n{\bf Definition.}  For $w \in \widehat{\R}$ and $\zeta >0$ define
   $$
   \Phi_\zeta(w) = \zeta^{1/3} \Phi(w) - \frac{1}{1+w^2}.
   $$
   
\begin{lemma} \label{lem Phi zeta}  Suppose $\mu > 0$ and let $\zeta =  \dfrac{b^2 \sigma^2}{2 a \mu^2}$.  We have 
   \begin{align*}
   \lefteqn{\left|  \int_M\left((2b\sigma)^{2/3} r ^{2/3} \Phi(w) - \frac{2a r^2}{1+w^2}\right)d\eta(r,w) - 2 \mu\int_M \Phi_\zeta (w) d\eta(r,w)\right|} \hspace{40ex} \\
   &  \le (2 b \sigma)^{2/3} \left(\frac{2\sigma^2}{a}\right)^{1/6}+2a \left(\frac{2\sigma^2}{a}\right)^{1/2}.
  \end{align*}
\end{lemma}  

\n{\bf Proof.}  Notice first
  $$
  (2b\sigma)^{2/3} r_0^{2/3} \Phi(w) - \frac{2a r_0^2}{1+w^2} = 2 \mu \zeta^{1/3}\Phi(w) -\frac{ 2 \mu}{1+w^2} = 2\mu \Phi_\zeta(w),
  $$  and so
\begin{align*}
  \left|(2b\sigma)^{2/3} r ^{2/3} \Phi(w) - \frac{2a r^2}{1+w^2} - 2 \mu \Phi_\zeta(w)\right|
   & = \left|(2b\sigma)^{2/3}\left(r^{2/3}-r_0^{2/3}\right) \Phi(w)- \frac{2a(r^2-r_0^2)}{1+w^2} \right|\\
  & \le (2b\sigma)^{2/3}\left|r^{2/3}-r_0^{2/3}\right|+ 2a\left|r^2-r_0^2\right|.
  \end{align*} 
Using Lemma \ref{lem nu bounds2} in the Appendix we have 
   \begin{equation} \label{r2 est}
   \int_{\R_+}\left|r^2-r_0^2\right|d\nu_{(\mu,a,\sigma)}(r) \le \left(\int_{\R_+} \left(r^2-r_0^2\right)^2 d\nu_{(\mu,a,\sigma)}(r)\right)^{1/2} \le \left(\frac{2\sigma^2}{a}\right)^{1/2}
  \end{equation} and similarly
 \begin{equation} \label{r23 est}
   \int_{\R_+}\left|r^{2/3}-r_0^{2/3}\right|d\nu_{(\mu,a,\sigma)}(r) \le  \int_{\R_+}\left|r^2-r_0^2\right|^{1/3}d\nu_{(\mu,a,\sigma)}(r) %\le\left(\int_{\R_+} \left(r^2-r_0^2\right)^2 d\nu_{(\mu,a,\sigma)}(r)\right)^{1/2}
    \le \left(\frac{2\sigma^2}{a}\right)^{1/6}
  \end{equation} 
and the result follows directly.   \qed

\s

\begin{lemma} \label{lem M CE} Suppose $\mu > 0$ and let $\zeta =  \dfrac{b^2 \sigma^2}{2 a \mu^2}$.  The generator ${\cal M}$ for the diffusion $\{(r_t,w_t): t \ge 0\}$ can be written in the form
 \begin{equation} \label{M CE}
 {\cal M} = 2 \mu {\cal M}_{3,\zeta} + (2b\sigma)^{2/3}\left(r^{2/3} - r_0^{2/3}\right) {\cal M}_4 + a(r^2-r_0^2){\cal M}_5 + \sigma^2 {\cal M}_6 + \frac{\sigma^4}{(2 b \sigma)^{2/3}}
 {\cal M}_7
\end{equation}
where 
     \begin{align}
   {\cal M}_{3,\zeta}  & = \zeta^{1/3}\left((1+w^3)\frac{\partial}{\partial w} + \frac{w^4}{2}\frac{\partial^2}{\partial w^2}\right) +w \frac{\partial}{\partial w} \nonumber \\
   {\cal M}_4 & = (1+w^3)\frac{\partial}{\partial w} + \frac{w^4}{2}\frac{\partial^2}{\partial w^2} \nonumber\\
    {\cal M}_5 & = -r \frac{\partial}{\partial r} + \frac{10 w}{3} \frac{\partial}{\partial w} \nonumber \\
    {\cal M}_6 & = \frac{1}{2r}\frac{\partial}{\partial r} + \frac{17 w}{9r^2} \frac{\partial}{\partial w} + \frac{1}{2}\frac{\partial^2}{\partial r^2}  - \frac{4w}{3r}\frac{\partial^2}{\partial r \partial w} + \frac{17 w^2}{9r^2}\frac{\partial^2}{\partial w^2} \nonumber\\
    {\cal M}_7& = \frac{1}{2r^{14/3}}\frac{\partial^2}{\partial w^2}. \nonumber
   \end{align}
\end{lemma}

\n{\bf Proof.}  This is simply a rearrangement of the expression \eqref{M} for ${\cal M}$ in Corollary \ref{cor lam}.  \qed

\s

Next we adapt and extend the arguments of Lemma \ref{lem G} and Proposition \ref{prop Phi}.  Recall the function $\Psi(\zeta)$ given in \eqref{Psi}.

\begin{lemma} \label{lem G CE}  For each $\zeta > 0$ there is a smooth bounded function $G_\zeta(w)$ such that
       \begin{equation} \label{G CE}
    {\cal M}_{3,\zeta}G_\zeta(w) =    \Phi_\zeta(w) - \Psi(\zeta).
    \end{equation}
Moreover for any compact set $D_1 \subset (0,\infty)$ the $G_\zeta(w)$ may be chosen so that $G_\zeta(w)$ and $(1+w^2) G_\zeta'(w)$ and $(1+w^2) G_\zeta''(w)$ are bounded uniformly for $\zeta \in D_1$.
 \end{lemma} 
 
 \n{\bf Proof.}  Consider the linear SDE \eqref{YB}:
  $$
            dY_t = \begin{bmatrix} -1  & 0 \\ \zeta^{1/3} & 0 \end{bmatrix}Y_t dt + \begin{bmatrix} 0 & \zeta^{1/6} \\ 0 & 0 \end{bmatrix} Y_t dW_t. 
     $$ 
Writing $Y_t = \|Y_t\| \begin{bmatrix} \cos \vartheta_t \\ \sin \vartheta_t \end{bmatrix}$ and then $w_t = \tan \vartheta_t$ we have 
    \begin{equation} \label{dwB}
    dw_t = \Bigl(\zeta^{1/3}(1+w_t^3) + w_t\Bigr)dt - \zeta^{1/6}w_t^2 dW_t
   \end{equation}
and 
   $$
   d\log \|X_t\| = \left(\zeta^{1/3} \Phi(w_t) - \frac{1}{1+w_t^2}\right)dt + \zeta^{1/6}\frac{w_t}{1+w_t^2}dW_t.
   $$  
Thus ${\cal M}_{3,\zeta}$ is the generator for the diffusion $\{w_t: t \ge 0\}$ given by \eqref{dwB}, and $\Phi_\zeta(w)$ is the integrand in the Furstenberg-Khasminskii formula for the top Lyapunov exponent $\Psi(\zeta)$ of \eqref{YB}.  Therefore ${\cal M}_{3,\zeta}G_\zeta(w) = \Phi_\zeta(w) - \Psi(\zeta)$ is the adjoint equation for \eqref{YB}.  By Arnold, Oeljeklaus and Pardoux \cite[Thm 5.1]{AOP86} we have the existence of $G_\zeta$, and by Imkeller and Lederer \cite[Theorem 3]{IL01} we have the formula \eqref{Psi} for $\Psi(\zeta)$.  Moreover there is a smooth $\pi$-periodic function $\widetilde{G}_\zeta$ such that $\widetilde{G}_\zeta(\vartheta) = G_\zeta(\tan \vartheta)$.  By Baxendale \cite[Prop.\ 4.2]{Bax94} the $\widetilde{G}_\zeta$ may be chosen so that the mapping $\zeta \mapsto \widetilde{G}_\zeta$ of $(0,\infty)$ to $C^2(\R/\pi \Z)$ is continuous.  In particular the $C_2$ norm of $\widetilde{G}_\zeta$ is uniformly bounded for $\zeta \in D_1$, and the proof is completed as in Lemma \ref{lem G}.  \qed

\begin{proposition}  \label{prop Phi CE}    Fix $\mu > 0$ and $a >0$ and a compact set $D \subset (0,\infty)$.  For all $\e > 0$ there exists $\sigma_0 > 0$ such that 
    \begin{equation} \label{Phi CE}
 \left|\int_M \Phi_\zeta(w)d\eta(r,w)  - \Psi(\zeta) \right| <\e
     \end{equation}
whenever $b\sigma \in D$ and $0 < \sigma \le \sigma_0$, where $\zeta =  \dfrac{b^2 \sigma^2}{2 a \mu^2}$. 
   \end{proposition}
       
\n{\bf Proof.} In the expression \eqref{M CE} for ${\cal M}$ some of the coefficients in ${\cal M}_6$ and ${\cal M}_7$ are not integrable functions of $r$.  As in the proof of Proposition \ref{prop Phi} we use truncation.  Let $C \subset (0,\infty)$ be a neighborhood of $r_0$ and let $h:(0,\infty) \to [0,1]$ be a smooth compactly supported function with $h\big|_{C} = 1$.  Define $G_{\zeta,h}(r,w) = h(r)G_\zeta(w)$. Then Lemmas \ref{lem M CE} and \ref{lem G CE} give
  \begin{align*}
   {\cal M} G_{\zeta,h}(r,w) &  =  2 \mu h(r) \bigl(\Phi_\zeta(w)-\Psi(\zeta)\bigr) + (2b\sigma)^{2/3}\bigl(r^{2/3} - r_0^{2/3}\bigr) {\cal M}_4G_{\zeta,h}(r,w)\\
   & \qquad + a(r^2-r_0^2){\cal M}_5G_{\zeta,h}(r,w) + \sigma^2 {\cal M}_6G_{\zeta,h}(r,w) + \frac{\sigma^4}{(2 b \sigma)^{2/3}}
 {\cal M}_7G_{\zeta,h}(r,w).
   \end{align*} 
Let $D_1 = \{\zeta= b^2\sigma^2/(2a \mu^2): b\sigma \in D\}$.  By Lemma \ref{lem G CE} the functions $G_\zeta(w)$, $wG_\zeta'(w)$, $G_\zeta''(w)$ and $w^2G_\zeta''(w)$ are bounded uniformly for $\zeta \in D_1$.  Since ${\cal M}_4=\zeta^{-1/3}\left({\cal M}_3 - w \frac{\partial}{\partial w}\right)$ then ${\cal M}_4 G_{\zeta,h}(r,w) = \zeta^{-1/3}h(r)\left(\Phi_\zeta(w) - \Psi(\zeta) - w G'_{\zeta}(w)\right)$ is bounded uniformly for $\zeta \in D_1$.  The terms ${\cal M}_5G_{\zeta,h}(r,w)$ and ${\cal M}_6G_{\zeta,h}(r,w)$ and ${\cal M}_7G_{\zeta,h}(r,w)$ involve a finite sum of products of the form $cf(r)g(w)$ where each $c$ is an absolute constant, and each $f(r)$ is a power of $r$ times a multiple of $h(r)$ or $h'(r)$ or $h''(r)$ and hence bounded, and each $g(w)$ is $G_\zeta(w)$ or  $wG_\zeta'(w)$ or $G_\zeta''(w)$ or $w^2G_\zeta''(w)$ and hence bounded uniformly for $\zeta \in D_1$.  Therefore ${\cal M}_5G_{\zeta,h}(r,w)$ and ${\cal M}_6G_{\zeta,h}(r,w)$ and ${\cal M}_7G_{\zeta,h}(r,w)$ are bounded in $(r,w)$ uniformly for $\zeta \in D_1$. Together there exists $K$ such that
        $
        |{\cal M}_j G_{\zeta,h}(r,w)| \le K
        $
for all $\zeta \in D_1$, $(r,w) \in M$, and $j=4,5,6,7$.  
Moreover since $G_\zeta$ and $\Phi_\zeta$ are bounded functions on $M$ it follows that $G_{\zeta,h}$ and ${\cal M}G_{\zeta,h}$ are bounded functions on $M$ and we may apply Proposition \ref{prop-bg} from the Appendix.  Thus  
   \begin{align*}
     0 & =  \int_M {\cal M}G_{\zeta,h}(r,w)d\eta(r,w) \\
    & =  2 \mu \int_M h(r) \bigl(\Phi_\zeta(w)-\Psi(\zeta)\bigr)d\eta(r,w) + (2b\sigma)^{2/3}\int_M \bigl(r^{2/3} - r_0^{2/3}\bigr) {\cal M}_4G_{\zeta,h}(r,w)d\eta(r,w)\\
   & \qquad + a\int_M (r^2-r_0^2){\cal M}_5G_{\zeta,h}(r,w)d\eta(r,w) + \sigma^2 \int_M {\cal M}_6G_{\zeta,h}(r,w) d\eta(r,w) \\
   & \qquad + \frac{\sigma^4}{(2 b \sigma)^{2/3}}
 \int_M {\cal M}_7G_{\zeta,h}(r,w)d\eta(r,w). 
     \end{align*} 
Now $|\Phi_\zeta(w) - \Psi(\zeta)| \le 2 \max\{|\Phi_\zeta(\tilde{w})|: \tilde{w} \in \widehat{\R}\} \le 2(\zeta^{1/3}+1)$.   Let $0 < d_1 < d_2$ be lower and upper bounds on the set $D$ and let $d_3$ be the corresponding upper bound on $D_1$.  Then if $b\sigma \in D$ we have
   \begin{align}
  \lefteqn{ 2 \mu  \left|\int_M \Phi_\zeta(w) d\eta(r,w) - \Psi(\zeta) \right| } \hspace{5ex} \nonumber \\
  % & = 2 \mu \left|\int_M \bigl( \Phi_\zeta(w) -\Psi(\zeta)\bigr)d\eta(r,w)  \right| \nonumber\\
   & \le  2 \mu \left|\int_M(1-h(r))\bigl( \Phi_\zeta(w) -\Psi(\zeta)\bigr)d\eta(r,w)  \right| + 2 \mu \left|\int_M  h(r) \bigl( \Phi_\zeta(w) -\Psi(\zeta)\bigr)d\eta(r,w)  \right|    \nonumber \\
   & \le  4 \mu(\zeta^{1/3}+1) \nu_{(\mu,a,\sigma)}(\R_+ \setminus C) + (2b\sigma)^{2/3} \left|\int_M \bigl(r^{2/3} - r_0^{2/3}\bigr) {\cal M}_4G_{\zeta,h}(r,w)d\eta(r,w)\right| \nonumber \\
   & \qquad + a \left|\int_M (r^2-r_0^2){\cal M}_5G_{\zeta,h}(r,w)d\eta(r,w)\right| + \sigma^2 \left|\int_M{\cal M}_6G_{\zeta,h}(r,w) d\eta(r,w)\right| \nonumber  \\
   & \qquad + \frac{\sigma^4}{(2 b \sigma)^{2/3}}
 \left|\int_M {\cal M}_7G_{\zeta,h}(r,w)d\eta(r,w)\right| \nonumber \\
    & \le  4 \mu(d_3^{1/3}+1) \nu_{(\mu,a,\sigma)}(\R_+ \setminus C) + (2d_2)^{2/3}K\int_{\R_+} \big|r^{2/3} - r_0^{2/3}\big| d\nu_{(\mu,a,\sigma)}(r) \nonumber \\
   & \qquad + a K \int_{\R_+} \big|r^2-r_0^2 \big|d\nu_{(\mu,a,\sigma)}(r) + \sigma^2 K  + \frac{\sigma^4}{(2 d_1)^{2/3}}K.  \label{Phi CE est}
 \end{align}
For fixed $\mu > 0$ and $a > 0$ and $D$, the right side of \eqref{Phi CE est} depends only on $\sigma$.  Using \eqref{r23 est} and \eqref{r2 est} on the second and third terms, and Markov's inequality on the first term, we see that the right side of \eqref{Phi CE est} converges to 0 as $\sigma \to 0$, and the proof is complete.   \qed

 \s
 
 \n{\bf Proof of Theorem \ref{thm CE unif}.}  Fix $\mu > 0$ and $a > 0$ and a compact set $D \subset (0,\infty)$ and write $\zeta = b^2\mu^2/(2\mu^2 a)$.  Since $\mu$ and $a$ are fixed we abbreviate $\nu_{(\mu,a,\sigma)} = \nu_\sigma$ to emphasize the dependence of the invariant measure $\nu$ on the parameter $\sigma$.  %Recall $r_0 = \sqrt{\mu/a}$. 
  By Corollary \ref{cor lam} and Lemma \ref{lem Phi zeta}
  % in the Appendix and (\ref{int r23}, \ref{r2 est}), it is enough to show that for all $\e > 0$ there exists $\sigma_0 > 0$ such that 
  we have 
   \begin{align*}
   \lefteqn{\left|\lambda(\mu,a,b,\sigma) - 2 \mu \Psi\left(\frac{b^2\sigma^2}{2 \mu^2 a}\right)\right|} \hspace{5ex} \\
   & \le \left|\lambda(\mu,a,b,\sigma) - (2b\sigma)^{2/3}\int_M r ^{2/3} \Phi(w)d\eta(r,w)+ 2a\int_M \frac{r^2}{1+w^2}d\eta(r,w)\right| \\
   & \quad +\left|  \int_M\left((2b\sigma)^{2/3} r ^{2/3} \Phi(w) - \frac{2a r^2}{1+w^2}\right)d\eta(r,w) - 2 \mu\int_M \Phi_\zeta (w) d\eta(r,w)\right|\\
   & \quad + 2\mu  \left|\int_M \Phi_\zeta(w) d\eta(r,w) - \Psi(\zeta) \right|\\
     & \le \frac{59 \sigma (2b\sigma)^{1/3}}{36} \int_{\R_+} r^{-2/3}d\nu_\sigma(r) + \frac{a  \sigma}{(2b\sigma)^{1/3}}\int_{\R_+} r^{2/3}d\nu_\sigma(r) +\frac{7}{3}\int_{\R_+}|\mu-ar^2|d\nu_\sigma(r)\\
     & \quad + (2 b \sigma)^{2/3} \left(\frac{2\sigma^2}{a}\right)^{1/6}+2a \left(\frac{2\sigma^2}{a}\right)^{1/2}  + 2\mu  \left|\int_M \Phi_\zeta(w) d\eta(r,w) - \Psi(\zeta) \right|.
   \end{align*}
The result now follows by using Lemma \ref{lem nu bounds} on the first and second terms, the inequality \eqref{r2 est} on the third term, and Proposition \ref{prop Phi CE} on the final term.
  \qed

\section{Proof of Theorem \ref{thm sigma0}} \label{sec sigma0 proof}

The proof uses the adjoint method in the more usual setting of an asymptotic expansion.  We work with the original variables $(r,\psi) \in M = (0,\infty) \times \R/(\pi \Z)$ of Section \ref{sec eval} and do not use the tangent transform of Section \ref{sec tangent}. 

 Write 
     $
     {\cal L} = {\cal L}_0+\sigma^2 {\cal L}_2
     $
      where 
     $$
      {\cal L}_0 = \left(\mu r- ar^3\right)\frac{\partial}{\partial r} + 2r^2 \cos \psi\bigl(b \cos \psi + a \sin \psi\bigr) \frac{\partial}{\partial \psi}
    $$
and 
$$
   {\cal L}_2 = \frac{1}{2r} \frac{\partial}{\partial r}+  \frac{1}{2}\frac{\partial^2}{\partial r^2} + \frac{1 }{2r^2}\frac{\partial^2}{\partial \psi^2},
  $$
and recall $Q(r,\psi) = \mu-ar^2 +2r^2 \cos \psi \bigl( b \sin \psi- a \cos  \psi \bigr)$.  Taking $F(r,\psi) = F_0(r,\psi)+\sigma^2F_2(r,\psi) +\cdots$ and $\lambda = \lambda_0 + \sigma^2 \lambda_2+ \cdots$ the adjoint equation becomes 
  $$({\cal L}_0+\sigma^2{\cal L}_2)(F_0+\sigma^2 F_2+\cdots)(r,\psi) = Q(r,\psi) - (\lambda_0+\sigma^2 \lambda_2+ \cdots)
  $$
 so that
   \begin{align} \label{adj}
          {\cal L}_0F_0(r,\psi) & = Q(r,\psi) -\lambda_0,  \\
    {\cal L}_0F_2(r,\psi)+{\cal L}_2 F_0(r,\psi) & = - \lambda_2, \label{adj2}
    \end{align}
   and so on.  

The fact that the operator ${\cal L}_0$ is not elliptic makes the existence of $F_0$ solving \eqref{adj} and then the existence of $F_2$ solving \eqref{adj2} more complicated than is typical.  Since ${\cal L}_0$ is a first order differential operator it may be interpreted as a vector field $V(r,\psi) = \begin{bmatrix} \mu r- ar^3 \\ 2r^2 \cos \psi(b \cos \psi + a \sin \psi)\end{bmatrix}$ on $M = (0,\infty) \times \R/\pi \Z$. Let $(r,\psi) \mapsto \Xi_t(r,\psi)$ denote the time $t$ flow on $M$ along the vector field $V$.  Write $r_0 = \sqrt{\mu/a}$ and $\psi_\ast = \arctan(-b/a)$.  The flow $\{\Xi_t: t \ge 0\}$ has a stable fixed point at $(r_0,\pi/2)$ and a hyperbolic fixed point at $(r_0,\psi_\ast)$.  For $0 < r < \infty$ and $\psi \neq \psi_\ast$ the path $\Xi_t(r,\psi)$ converges to $(r_0,\pi/2)$ as $t \to \infty$.  Moreover $DV(r_0,\pi/2) = \begin{bmatrix} -2\mu & 0 \\ 0 & -2\mu \end{bmatrix}$ so that the convergence is exponentially fast for each point $(r,\psi)$ with $\psi \neq \psi_\ast$. 

Consider the equation \eqref{adj}.  Choose 
 \begin{equation} \label{lambda0}
  \lambda_0 = Q(r_0,\pi/2) = \mu-ar_0^2 = 0.
  \end{equation} 
  Then $Q(\Xi_t(r,\psi)) \to Q(r_0,\pi/2) = \lambda_0 $ exponentially quickly as $t \to \infty$, and we get the solution 
     \begin{equation} \label{F0}
     F_0(r,\psi) = -\int_0^\infty \Bigl(Q(\Xi_t(r,\psi))-\lambda_0 \Bigr)dt 
     \end{equation}  
valid for $0 < r < \infty$ and $\psi \neq \psi_\ast$.  

Next we choose $\lambda_2$ so that \eqref{adj2} has a solution $F_2$.  Arguing as above we require 
   $$
   \lambda_2 = -{\cal L}_2F_0(r_0,\pi/2).
   $$ 
The evaluation of ${\cal L}_2F_0(r_0,\pi/2)$ involves an elementary local analysis.  The function $f(r) = F_0(r,\pi/2)$ satisfies
     $$
     (\mu r - ar^3)f'(r) = {\cal L}_0F_0(r,\pi/2) = Q(r,\pi/2)-\lambda_0 = \mu - ar^2.
     $$
Thus $f'(r) = 1/r$ and $f''(r) = -1/r^2$.  %In particular $f''(r_0) = -1/r_0^2 = - a/\mu$.  
Similarly the function $g(\psi) = F_0(r_0,\psi)$ satisfies
        $$
        2r_0^2\cos \psi \bigl(b \cos \psi+a \sin \psi\bigr)g'(\psi) = {\cal L}_0F_0(r_0,\psi) = Q(r_0,\psi) -\lambda_0 = 2r_0^2 \cos \psi\bigl(b\sin \psi-a \cos \psi\bigr).
        $$  Thus   
   $$
     g'(\psi) = \frac{ b\sin \psi-a \cos \psi}{b \cos \psi+a \sin \psi}
     $$   
 and so 
  $$
   g''(\psi) = \frac{ a^2+b^2}{(b \cos \psi+a \sin \psi)^2}.
   $$       
Then
  \begin{align}
  \lambda_2 = -{\cal L}_2 F_0(r_0,\psi_0) & =  -\left(\frac{1}{2r_0} \frac{\partial F_0}{\partial r}(r_0,\psi_0) + \frac{1}{2} \frac{\partial^2 F_0}{\partial r^2}(r_0,\psi_0) + \frac{1}{2r_0^2}  \frac{\partial^2 F_0}{\partial \psi^2}(r_0,\psi_0)\right) \nonumber\\
  &  = -\left( \frac{1}{2r_0}f'(r_0) + \frac{1}{2}f''(r_0)+ \frac{1}{2r_0^2}g''(\psi_0) \right)\nonumber \\
  & = -\left( \frac{1}{2r_0}\cdot\frac{1}{r_0} + \frac{1}{2}\cdot\frac{-1}{r_0^2}  + \frac{1}{2r_0^2} \cdot \frac{a^2+b^2}{a^2}\right) \nonumber \\
  & =  -\frac{a^2+b^2}{2 \mu a}. \label{lambda2}
  \end{align}           
With this value of $\lambda_2$ the equation ${\cal L}_0F_2+{\cal L}_2 F_0 = - \lambda_2$ has the solution  
    \begin{equation} \label{F2}
     F_2(r,\psi) = \int_0^\infty \Bigl({\cal L}_2F_0(\Xi_t(r,\psi))+\lambda_2 \Bigr)dt 
     \end{equation}
valid for $0 < r < \infty$ and $\psi \neq \psi_\ast$.

\begin{remark}  The argument above describes a general method for finding the functions $F_0$ and $F_2$.  In this case there is an explicit solution $F_0(r,\psi) = \log r - \log|b\cos \psi+a \sin \psi|$ for $0 < r < \infty$ and $\psi \neq \psi_\ast$.  Then 
    $${\cal L}_2F_0(r,\psi) = \frac{a^2+b^2}{2r^2(b \cos \psi+ a \sin \psi)^2},
    $$ 
and there is a more direct calculation 
   $$\lambda_2= -{\cal L}_2F_0(r_0,\pi/2) = -\frac{a^2+b^2}{2r_0^2a^2} = -\frac{a^2+b^2}{2\mu a}.
    $$ \end{remark}

%The adjoint method is based on the idea that ${\cal L}F(r,\psi) = Q(r,\psi)-\lambda$ implies
%    $$
%    0 = \int_M {\cal L}F(r,\psi) d\widetilde{\nu}(r,\psi) = \int_M Q(r,\psi) d\widetilde{\nu}(r,\psi) - \lambda.
%    $$
Replacing $F$ by $F_0+\sigma^2 F_2$ in the adjoint equation and using (\ref{adj},\ref{adj2}) gives
    \begin{align*}
    {\cal L}(F_0+\sigma^2 F_2)(r,\psi) & = ({\cal L}_0 + \sigma^2 {\cal L}_2)(F_0+\sigma^2 F_2)(r,\psi) \\
    & = {\cal L}_0 F_0(r,\psi)+\sigma^2({\cal L}_0F_2+{\cal L}_2 F_0)(r,\psi)+\sigma^4 {\cal L}_2F_2(r,\psi)\\
    & = Q(r,\psi) - \lambda_0 -\sigma^2 \lambda_2 + \sigma^4 {\cal L}_2F_2(r,\psi)
    \end{align*}
and integrating with respect to $\widetilde{\nu}$ should give
   $$
    \left|\int_M Q(r,\psi) d\widetilde{\nu}(r,\psi) - \lambda_0-\sigma^2 \lambda_2 \right| \le  \sigma^4 \int_M\left|{\cal L}_2 F_2(r,\psi) \right|d\widetilde{\nu}(r,\psi).
    $$
However in the present setting the functions $F_0$ and $F_2$ are undefined at $\psi_\ast$ and may be unbounded as $r \to 0$ and as $r \to \infty$.  There is no guarantee that $\int_M {\cal L}(F_0+\sigma^2 F_2)d\widetilde{\nu} = 0$, and no guarantee that $\int_M |{\cal L}_2F_2|d\widetilde{\nu}$ is bounded uniformly as $\sigma \to 0$.  The proof that follows replaces $F_0$ and $F_2$ by smooth functions with compact support on $(0,\infty) \times (\R/\pi\Z \setminus\{\psi_\ast\})$, and uses an {\it a priori} estimate on the invariant measure $\widetilde{\nu}$ to quantify the effect of these replacements.

The next two lemmas deal with the dependence of the functions $F_0$ and $F_2$ on the parameter $b$, and with the dependence of the invariant probability $\widetilde{\nu}= \widetilde{\nu}_{(b,\sigma)}$ on the parameters $\sigma$ and $b$.  (Recall that the other two parameters $\mu > 0$ and $a >0$ are fixed.)  For ease of notation we do not explicitly write the dependence of the unstable fixed point $\psi_\ast = \arctan(-b/a)$ and the functions $F_0$ and $F_2$ and the set $\widehat{U}$ on the parameter $b$. The proofs of the two lemmas will be given after the proof of Theorem \ref{thm sigma0}.

For $\psi, \psi' \in \R/{\pi \Z}$ let $d(\psi,\psi')$ denote the geodesic distance $\min\{|\psi-(\psi'+n\pi)|: n \in \Z\}$.

\begin{lemma}  \label{lem unif} Fix $\mu > 0$ and $a> 0$ and $k < \infty$.  For all $0 < r_1 < r_2 <\infty$ and $\e > 0$ there exists $K <\infty$ such that 
  $$
  \sup_{\substack{r_1 \le r \le r_2,\\d(\psi,\psi_\ast) \ge \e}} \max\left\{\bigl|F_i(r,\psi)\bigr|,
  \left|\frac{\partial F_i}{\partial r}(r,\psi)\right|,\left|\frac{\partial^2 F_i}{\partial r^2}(r,\psi)\right|,
  \left|\frac{\partial F_i}{\partial \psi}(r,\psi)\right|,\left|\frac{\partial^2 F_i}{\partial \psi^2}(r,\psi)\right| \right\} \le K
  $$
for $i = 0,2$ whenever $0 \le b \le k$.
\end{lemma}

\begin{lemma} \label{lem nu U}  Fix $\mu > 0$ and $a > 0$ and $k  > 0$.  For all $0 < s_1 < \sqrt{\mu/a} < s_2 < \infty$ and $\sigma_0 > 0$ there exist constants $c_1< \infty$, $c_2>0$ and $\delta > 0$ such that  
     $$
     \int_{M \setminus M_0} (1+r^2) d\widetilde{\nu}_{b,\sigma}(r,\psi) \le c_1e^{-c_2/\sigma^2}
     $$
whenever $0 < \sigma \le \sigma_0$ and $0 \le b \le k$, where 
   $$
    M_0 = [s_1,s_2] \times \{\psi \in \R/\pi\Z: d(\psi,\psi_\ast) \ge \delta \}.
 $$
\end{lemma}
         
\n{\bf Proof of Theorem \ref{thm sigma0}.}   Fix $0 < r_1 < s_1 < \sqrt{\mu/a} < s_2 < r_2 < \infty$, and let $c_1$, $c_2$ and $\delta$ be the constants given by Lemma \ref{lem nu U}.  Choose smooth $h:(0,\infty) \to [0,1]$ such that $h(r) = 0$ for $r \not\in[r_1,r_2]$ and $h(r) =1$ for $r \in [s_1,s_2]$.  Choose smooth $j: \R/\pi\Z \to [0,1]$ such that $j(\psi)= 0$ for $d(\psi,0) < \delta/2$ and $j(\psi) = 1$ for $d(\psi,0) \ge \delta$.  Then define 
    $$
    \widehat{F_i}(r,\psi) = h(r)j(\psi-\psi_\ast) F_i(r,\psi)
    $$
for $i = 0,2$.  
On the set $M_0$ we have 
   \begin{align*}
    {\cal L}_0\widehat{F}_0 & = {\cal L}_0 F_0 = Q-\lambda_0, \\
    {\cal L}_0\widehat{F}_2+{\cal L}_2\widehat{F}_0 & = {\cal L}_0F_2+{\cal L}_2 F_0 = -\lambda_1.
 \end{align*}
Also the functions ${\cal L}_0 \widehat{F}_0$ and ${\cal L}_0\widehat{F}_2+{\cal L}_2\widehat{F}_0$ and ${\cal L}_2\widehat{F}_2$ all have support in the set $[r_1,r_2] \times \{\psi \in \R/\pi\Z: d(\psi,\psi_\ast) \ge \delta/2\}$.  Applying Lemma \ref{lem unif} with $\e = \delta/2$, we obtain constants $K_1$, $K_2$ and $K_3$ so that 
   \begin{equation} \label{K123}
   \begin{split}
    |{\cal L}_0\widehat{F}_0(r,\psi) -Q(r,\psi)+\lambda_0|  & \le K_1(1+r^2)\\
     |{\cal L}_0\widehat{F}_2(r,\psi)+{\cal L}_2\widehat{F}_0(r,\psi) +\lambda_2| & \le K_2 \\
    |{\cal L}_2\widehat{F}_2(r,\psi) | & \le K_3 
    \end{split}
   \end{equation}
for all $(r,\psi) \in M$ and all $0 \le b \le k$.  (The constants $K_1$, $K_2$ and $K_3$ depend on bounds on the first 2 derivatives of $h$ and $j$ as well as the constant $K$ from Lemma \ref{lem unif}.)  
 Now 
  $$
  {\cal L}(\widehat{F}_0+\sigma^2 \widehat{F}_2)  =  {\cal L}_0\widehat{F}_0 + \sigma^2({\cal L}_0\widehat{F}_2+{\cal L}_2\widehat{F}_0) + \sigma^4 {\cal L}_2\widehat{F}_2.
  $$
 Since $\widehat{F}_0+\sigma^2 \widehat{F}_2$ is a smooth function on $M$ with compact support, we may integate with respect to $\widetilde{\nu}_{(b,\sigma)}$ to obtain 
  \begin{align*}
  0  & =  \int_M{\cal L}(\widehat{F}_0+\sigma^2 \widehat{F}_2) d\widetilde{\nu}_{(b,\sigma)} \\
   & =\int_M{\cal L}_0\widehat{F}_0d\widetilde{\nu}_{(b,\sigma)} + \sigma^2\int_M({\cal L}_0\widehat{F}_2+{\cal L}_2\widehat{F}_0)d\widetilde{\nu}_{(b,\sigma)} + \sigma^4 \int_M{\cal L}_2\widehat{F}_2 d\widetilde{\nu}_{(b,\sigma)}\\
   & =  \int_M(Q-\lambda_0)d\widetilde{\nu}_{(b,\sigma)} + \sigma^2\int_M(-\lambda_2)d\widetilde{\nu}_{(b,\sigma)} + \sigma^4 \int_M{\cal L}_2\widehat{F}_2d\widetilde{\nu}_{(b,\sigma)}\\
     &  \quad  + \int_{M \setminus M_0}({\cal L}_0\widehat{F}_0-Q+\lambda_0)d\widetilde{\nu}_{(b,\sigma)} + \sigma^2\int_{M \setminus M_0}({\cal L}_0\widehat{F}_2+{\cal L}_2\widehat{F}_0+\lambda_2)d\widetilde{\nu}_{(b,\sigma)} \\
  & =  \lambda(\mu,a,b,\sigma) +  \frac{\sigma^2(a^2+b^2)}{2 \mu a} + \sigma^4 \int_M{\cal L}_2\widehat{F}_2d\widetilde{\nu}_{(b,\sigma)}\\
     &  \quad  + \int_{M \setminus M_0}({\cal L}_0\widehat{F}_0-Q+\lambda_0)d\widetilde{\nu}_{(b,\sigma)} + \sigma^2\int_{M \setminus M_0}({\cal L}_0\widehat{F}_2+{\cal L}_2\widehat{F}_0+\lambda_2)d\widetilde{\nu}_{(b,\sigma)},
 \end{align*}
where the last line uses the Furstenberg-Khasminskii formula and (\ref{lambda0},\ref{lambda2}).  Rearranging and using \eqref{K123} and Lemma \ref{lem nu U} we get
   \begin{align*} 
    \left|\lambda(\mu,a,b,\sigma) +  \frac{\sigma^2(a^2+b^2)}{2 \mu a}\right| 
       & \le  \sigma^4 \int_M \left|{\cal L}_2\widehat{F}_2 \right|d\widetilde{\nu}_{(b,\sigma)}+ \int_{M \setminus M_0}\left|{\cal L}_0\widehat{F}_0-Q+\lambda_0\right|d\widetilde{\nu}_{(b,\sigma)}\\
     & \qquad  + \sigma^2\int_{M \setminus M_0}\left|{\cal L}_0\widehat{F}_2+{\cal L}_2\widehat{F}_0+\lambda_2\right|d\widetilde{\nu}_{(b,\sigma)} \\
     & \le  \sigma^4  K_3 + \int_{M \setminus M_0}[K_1(1+r^2)+\sigma^2 K_2]d\widetilde{\nu}_{(b,\sigma)}\\
     & \le \sigma^4 K_3 + (K_1+\sigma^2 K_2)c_1e^{-c_2/\sigma^2}. 
     \end{align*}
Using the estimates $xe^{-x} \le e^{-1}$ and $x^2e^{-x} \le 4e^{-2}$ with $x =c_2/\sigma^2$ we get
  $$
  \left| \lambda+ \frac{\sigma^2(a^2+b^2)}{2 \mu a}\right| 
  \le \left(K_3+c_1\left(\frac{4K_1}{c_2^2e^2} + \frac{K_2}{c_2e}\right)\right)\sigma^4
  $$
and the proof of Theorem \ref{thm sigma0} is complete.  \qed   

\s

 \n{\bf Proof of Lemma \ref{lem unif}.} For any function $F(r,\psi)$ write
    $$
   S(F)(\psi) = \sup_{r_1 \le r \le r_2} \max\left\{|F(r,\psi)|,
  \left|\frac{\partial F}{\partial r}(r,\psi)\right|,\left|\frac{\partial^2 F}{\partial r^2}(r,\psi)\right|,
  \left|\frac{\partial F}{\partial \psi}(r,\psi)\right|,\left|\frac{\partial^2 F}{\partial \psi^2}(r,\psi)\right| \right\}.
  $$
In this proof write $\psi_\ast^{(b)}$ and $F_i^{(b)}$ to denote explicitly the dependence of the unstable fixed point $\psi_\ast$ and the functions $F_i$ on the parameter $b$.  Since the operator ${\cal L}_0$, equivalently the vector field $V(r,\psi)$, and the function $Q(r,\psi)$ are smooth functions of $(r,\psi)$ and $b$, then the function $F_0^{(b)}(r,\psi)$ defined in \eqref{F0} is a smooth function of $(r,\psi)$ and $b$ away from the unstable fixed point $\psi_\ast^{(b)}$.      More precisely, for fixed $b_0$ define the sets $C_{b_0} = \{\psi \in \R/\pi \Z: d(\psi,\psi_\ast^{(b_0)}) \ge 2\e/3\}$ and $U_{b_0} = \{b \in \R : d(\psi_\ast^{(b)},\psi_\ast^{(b_0)}) < \e/3\}$.  Then $(r,\psi,b) \to F_0^{(b)}(r,\psi)$ is a smooth mapping of $[r_1,r_2]\times C_{b_0} \times \overline{U}_{b_0}$ into $\R$.  It follows that there is $K_{b_0} <\infty$ such that 
     $$
 \sup_{b \in U_{b_0}} \sup_{\psi \in C_{b_0}} S(F_0^{(b)})(\psi) %\max\left\{|F_0^{(b)}(r,\psi)|,
 % \left|\frac{\partial F_0^{(b)}}{\partial r}(r,\psi)\right|,\left|\frac{\partial^2 F_0^{(b)}}{\partial r^2}(r,\psi)\right|,
 % \left|\frac{\partial F_0^{(b)}}{\partial \psi}(r,\psi)\right|,\left|\frac{\partial^2 F_0^{(b)}}{\partial \psi^2}%(r,\psi)\right| \right\} 
 \le K_{b_0}
  $$
for all $b \in U_{b_0}$.  If $b \in U_{b_0}$ and $d(\psi,\psi_\ast^{(b)}) \ge \e$ then 
  $$
  d(\psi, \psi_\ast^{(b_0)}) \ge d(\psi,\psi_\ast^{(b)}) - d(\psi_\ast^{(b)},\psi_\ast^{(b_0)}) > \e - \e/3 = 2\e/3.
  $$  
It follows that 
   $$
   \sup_{b \in U_{b_0}} \sup_{d(\psi, \psi_\ast^{(b)}) \ge \e} S(F_0^{(b)})(\psi) \le 
        \sup_{b \in U_{b_0}} \sup_{\psi \in C_{b_0}} S(F_0^{(b)})(\psi) \le K_{b_0},
   $$     
and the result for $F_0^{(b)}$ now follows by compactness of the set $[0,k]$.  The result for $F_2^{(b)}$ is similar, using the smoothness in $(r,\psi)$ and $b$ of ${\cal L}_2 F_0^{(b)}$ in the formula \eqref{F2} for $F_2^{(b)}$.    
\qed 

\s

\n {\bf Proof of Lemma \ref{lem nu U}.}  Recall $\mu > 0$ and $a >0$ are fixed, and $0 < s_1< \sqrt{\mu/a} < s_2 <\infty$.  Write $\nu = \nu_\sigma$ to denote its dependence on the parameter $\sigma$.  Taking $\gamma = a/(4\sigma^2)$ in Lemma \ref{lem nu bounds2} in the Appendix gives
      $$
      \int_{\R_+} e^{a (r^2-\mu/a)^2/(4 \sigma^2)} d\nu_\sigma(r) \le 2\sqrt{2}.
      $$
It follows from Markov's inequality that there exist constants $c_3 < \infty$ and $c_4 >0$ such that 
    \begin{equation} \label{nu r ineq}
     \int_{r \not\in[s_1,s_2]} (1+r^2) d\nu_\sigma(r)   \le c_3e^{-c_4/\sigma^2}
    \end{equation} 
whenever $0 < \sigma \le \sigma_0$.

Let $f: \R/\pi \Z \to [0,\infty)$ be $C^2$ with $f(\psi) = d(\psi,0)^2$ for $d(\psi,0)\le 1$, and $f(\psi) \ge 1$ otherwise.  Given parameter values $b$ and $\sigma$, define the function $H_{b,\sigma}: \R/\pi \Z \to [0,1]$ by 
   $$
   H_{b,\sigma}(\psi)  = \exp\left( \frac{-\gamma f(\psi-\psi_\ast)}{\sigma^2}\right)
   $$
where $\gamma >0$ will be chosen later.  Let $h: \R_+ \to [0,1]$ be smooth with compact support such that $h(r) = 1$ for $s_1 \le r \le s_2$.  Write ${\cal L}_{b,\sigma}$ to emphasize the dependence of ${\cal L}$ on the parameters $b$ and $\sigma$.

\s

\n{\bf Claim:}  There exist positive constants $\gamma$, $\delta_1$, $A$ and $B$ such that 
  \begin{equation} \label{claim}
      \frac{{\cal L}_{b,\sigma}(hH_{b,\sigma})(r,\psi)}{H_{b,\sigma}(\psi)}  
   \le \begin{cases} -A & \mbox{ if } r \in [s_1,s_2] \mbox{ and } d(\psi,\psi_\ast) \le \delta_1 \\
                 B/\sigma^2 & \mbox{ otherwise}
                 \end{cases} 
      \end{equation}
whenever $0 < \sigma \le \sigma_0$ and $0 \le b \le k$. 

\s

\n{\bf Proof of Claim.}  We have
   \begin{align*}
  \lefteqn{ \frac{{\cal L}_{b,\sigma}(hH_{b,\sigma})(r,\psi)}{H_{b,\sigma}(r,\psi)}} \hspace{10ex} \\
   & =  \left(\mu r - ar^3+ \frac{\sigma^2}{2r}\right)h'(r)  + \frac{\sigma^2}{2}h''(r) \\
   & \quad +2r^2 h(r)\cos \psi\bigl(b\cos \psi+a\sin \psi\bigr) \frac{H_{b,\sigma}'(\psi)}{H_{b,\sigma}(\psi)}+\frac{\sigma^2 h(r)}{2r^2}\frac{H_{b,\sigma}''(\psi)}{H_{b,\sigma}(\psi)} \\
     & =  \left(\mu r - ar^3+ \frac{\sigma^2}{2r}\right)h'(r)  + \frac{\sigma^2}{2}h''(r)\\
     & \quad  + 2r^2 h(r) \cos \psi\bigl(b\cos \psi+a\sin \psi\bigr)\left(-\frac{\gamma}{\sigma^2} f'(\psi-\psi_\ast)\right)\\
     &  \quad + \frac{\sigma^2 h(r)}{2r^2}\left(-\frac{\gamma}{\sigma^2}f''(\psi-\psi_\ast)+\frac{\gamma^2}{\sigma^4}\bigl(f'(\psi-\psi_\ast)\bigr)^2\right)\\
   &  =   \left(\mu r - ar^3+ \frac{\sigma^2}{2r}\right)h'(r)  + \frac{\sigma^2}{2}h''(r)  -\frac{ \gamma h(r)}{2r^2}f''(\psi-\psi_\ast)\\
  &  \quad  + \frac{h(r)}{\sigma^2}\left(- 2\gamma r^2 \cos \psi\bigl(b\cos \psi+a\sin \psi\bigr) f'(\psi-\psi_\ast) +  \frac{ \gamma^2}{2r^2}(f'(\psi-\psi_\ast))^2\right).
 \end{align*}
The second line of \eqref{claim} follows immediately.   For $r \in [s_1,s_2]$ and $|\psi-\psi_\ast| \le 1$ we have 
  \begin{align*}
   \frac{{\cal L}_{b,\sigma}(hH_{b,\sigma})(r,\psi)}{H_{b,\sigma}(r,\psi)} 
   &  =    -\frac{ \gamma}{r^2}  + \frac{1}{\sigma^2}\left(- 4\gamma r^2 \cos \psi\bigl(b\cos \psi+a\sin \psi\bigr) (\psi-\psi_\ast) +  \frac{2 \gamma^2}{r^2}(\psi-\psi_\ast)^2\right).
 \end{align*}
The function $g(\psi) = \cos\psi(b \cos \psi+a \sin \psi)$ satisfies $g(\psi_\ast) = 0$ and $g'(\psi_\ast) = a$ and $|g''(\psi)| = |-2(b\sin 2\psi+a \sin 2 \psi)| \le 2\sqrt{a^2+b^2} \le 2 \sqrt{a^2+k^2}$.  Therefore for $|\psi-\psi_\ast| \le 1$ we have 
      $$
      |\cos\psi(b \cos \psi+a \sin \psi) - a(\psi-\psi_\ast)| \le \sqrt{a^2+k^2}(\psi-\psi_\ast)^2
      $$
and so
    $$
    \left|\cos \psi\bigl(b\cos \psi+a\sin \psi\bigr) (\psi-\psi_\ast) - a(\psi-\psi_\ast)^2\right| \le \sqrt{a^2+k^2}|\psi-\psi_\ast|^3.
    $$
Taking $\delta_1 = a/(2\sqrt{a^2+k^2}) < 1/2$ we get   
  $$
 \cos \psi\bigl(b\cos \psi+a\sin \psi\bigr) (\psi-\psi_\ast)  \ge \frac{a}{2}(\psi-\psi_\ast)^2
  $$
whenever $|\psi-\psi_\ast|  < \delta_1$.  Therefore 
  \begin{align*}
   \frac{{\cal L}_{b,\sigma}(hH_{b,\sigma})(r,\psi)}{H_{b,\sigma}(\psi)} 
   & \le   -\frac{ \gamma}{r^2} + \frac{1}{\sigma^2}\left(- 2\gamma r^2 a\,d(\psi,\psi_\ast)^2+  \frac{2 \gamma^2}{r^2}d(\psi,\psi_\ast)^2\right)\\
   & =   -\frac{ \gamma}{r^2} + \frac{2\gamma}{\sigma^2 r^2 }\bigl(\gamma - a r^4 \bigr)d(\psi,\psi_\ast)^2
\end{align*} 
whenever $r \in[s_1,s_2]$ and $d(\psi,\psi_\ast) \le \delta_1$. Choosing $0 < \gamma < a s_1^4$ gives
   $$
    \frac{{\cal L}_{b,\sigma}(hH_{b,\sigma})(r,\psi)}{H_{b,\sigma}(\psi)} 
    \le   -\frac{ \gamma}{s_2^2}
    $$
whenever $r \in [s_1,s_2]$ and $d(\psi,\psi_\ast) \le \delta_1$, and the proof of the claim is complete.

\s

Letting $U_1 = \{\psi \in \R/(\pi \Z): d(\psi,\psi_\ast) \le \delta_1\}$ we have 
   $$
   {\cal L}_{b,\sigma}(hH_{b,\sigma}) \le - A H_{b,\sigma}1_{[s_1,s_2]\times U_1} + \frac{B}{\sigma^2}H_{b,\sigma}1_{([s_1,s_2]\times U_1)^c}
   $$
The comparison theorem of Meyn and Tweedie \cite[Thm 1.1]{MTIII} gives
  $$
  A\int_0^t P_u\bigl(H_{b,\sigma}1_{[s_1,s_2]\times U_1}\bigr)(r,\psi)du \le h(r)H_{b,\sigma}(\psi) + \frac{B}{\sigma^2}\int_0^t P_u\bigl(H_{b,\sigma}1_{([s_1,s_2]\times U_1)^c}\bigr)(r,\psi)du
  $$
where $\{P_u: u \ge 0\}$ is the Markov semigroup corresponding to the generator ${\cal L}_{b,\sigma}$.  Since $h$ and $H_{b,\sigma}$ are bounded we can integrate with respect to the invariant probability measure $\widetilde{\nu}_{b,\sigma}$ and let $t \to \infty$ to obtain 
\begin{equation} \label{claim2}
   A\int_{[s_1,s_2] \times U_1}H_{b,\sigma}(\psi) d\widetilde{\nu}_{b,\sigma}(r,\psi) \le \frac{B}{\sigma^2}\int_{([s_1,s_2] \times U_1)^c } H_{b,\sigma}(\psi) d\widetilde{\nu}_{b,\sigma}(r,\psi).
      \end{equation}
The integral on the right side of \eqref{claim2} is 
   \begin{align*}
   \lefteqn{ \int_{([s_1,s_2] \times U_1)^c } H_{b,\sigma}(\psi) d\widetilde{\nu}_{b,\sigma}(r,\psi)} \hspace{12ex} \\
    & =  \int_{[s_1,s_2] \times U_1^c } H_{b,\sigma}(\psi) d\widetilde{\nu}_{b,\sigma}(r,\psi)
     + \int_{[s_1,s_2]^c\times \R/(\pi \Z) } H_{b,\sigma}(\psi) d\widetilde{\nu}_{b,\sigma}(r,\psi)\\
     & \le \int_{[s_1,s_2] \times U_1^c } e^{-\gamma \delta_1^2/\sigma^2} d\widetilde{\nu}_{b,\sigma}(r,\psi)
     + \int_{[s_1,s_2]^c \times \R/(\pi \Z)} 1 d\widetilde{\nu}_{b,\sigma}(r,\psi)\\
     & \le e^{-\gamma \delta_1^2/\sigma^2} + c_3 e^{-c_4/\sigma^2}
  \end{align*}   
where the last inequality uses \eqref{nu r ineq}.  Choose $0 < \delta <\min(\delta_1, \sqrt{c_2/\gamma})$ and define $U = \{\psi \in \R/\pi\Z: d(\psi,\psi_\ast)< \delta\} \subset U_1$.  For $\psi \in U$ we have $H_{b,\sigma}(\psi) \ge e^{-\gamma \delta^2/\sigma^2}$ and so 
 $$
   e^{-\gamma \delta^2/\sigma^2}\widetilde{\nu}_{b,\sigma}\bigl([s_1,s_2]\times U\bigr) 
     \le  \int_{[s_1,s_2]\times U} H_{b,\sigma}(r,\psi) d\widetilde{\nu}_{b,\sigma}(r,\psi)  \le  \frac{B}{A  \sigma^2 }\left(e^{-\gamma \delta_1^2/\sigma^2} + c_3e^{-c_4/\sigma^2}\right).
  $$
   Therefore 
   \begin{equation} \label{nu psi ineq}
   \widetilde{\nu}_{b,\sigma}\bigl([s_1,s_2]\times U\bigr) \le  \frac{B}{A \sigma^2 }\left(e^{-\gamma (\delta_1^2-\delta^2)/\sigma^2} + c_3e^{-(c_4-\gamma \delta^2)/\sigma^2}\right)
  \end{equation}
whenever $0 < \sigma \le \sigma_0$ and $0 \le b \le k$. Let $M_0 = [s_1,s_2] \times \{\psi \in \R/(\pi \Z): d(\psi,\psi_\ast) \ge \delta\} = [s_1,s_2]\times U^c$.  Then $M \setminus M_0 = \bigl([s_1,s_2]^c\times \R/(\pi \Z)\bigr) \cup \bigl([s_1,s_2]\times U\bigr)$, and 
%$\widehat{U} = ([s_1,s_2]^c \times \R/\pi\Z) \cup ([s_1,s_2] \times U)$ and 
using \eqref{nu r ineq} and \eqref{nu psi ineq} we get 
   \begin{align*}
   \int_{M\setminus M_0}(1+r^2)d\widetilde{\nu}_{b,\sigma}(r,\psi) 
       & =  \int_{r \not\in [s_1,s_2]}(1+r^2)d\nu_\sigma(r) + \int_{[s_1,s_2] \times U} (1+r^2)d\widetilde{\nu}_{b,\sigma}(r,\psi)\\
    & \le  c_3e^{-c_4/\sigma^2} + (1+s_2^2) \frac{B}{A \sigma^2 }\left(e^{-\gamma (\delta_1^2-\delta^2)/\sigma^2} + c_3e^{-(c_4-\gamma \delta^2)/\sigma^2}\right)
   \end{align*} 
whenever $0 < \sigma \le \sigma_0$ and $0 \le b \le k$.   Since $\dfrac{1}{\sigma^2}e^{-\kappa/\sigma^2} \le \dfrac{1}{\kappa e}$ for all positive $\kappa$ and $\sigma$, the result follows directly.  \qed

\appendix
\section*{Appendix}
\renewcommand{\thesubsection}{\Alph{subsection}}

%\section{Appendix} \label{sec-app}

\section{The stationary measure $\nu$} \label{sec nu}

The diffusion process $\{r_t: t \ge 0\}$ given by \eqref{r} has a unique invariant probability measure $\nu = \nu_{\mu,a,\sigma}$ on $\R_+ = (0,\infty)$ with density  
     \begin{equation} \label{rho}
  \rho(r) = \dfrac{2r \exp\left(-\frac{a}{2\sigma^2}\left[r^2-\frac{\mu}{a}\right]^2\right)}{\sqrt{\frac{\pi \sigma^2}{2a}}\mbox{erfc}\left(-\frac{\mu}{\sigma\sqrt{2a}}\right)}, \qquad 0 < r < \infty,
  \end{equation}
see \cite[eqn (21)]{DNR}.  Here
       $
      { \rm erfc}(x) = (2/\sqrt{\pi}) \int_x^\infty e^{-u^2}du
       $
is the complementary error function. %Note that ${\rm erfc}$ is monotone decreasing with $\lim_{x \to -\infty} {\rm erfc}(x) = 2$ and ${\rm erfc}(0) = 1$ and $\lim_{x \to \infty} {\rm erfc}(x) = 0$.  
It follows that $\int_{\R_+} r^\gamma d\nu(r) < \infty$ for all $\gamma > -2$.  There is an exact formula
  \begin{equation} \label{int r2}
   \int_{\R_+} r^2 d\nu(r) = \frac{\mu}{a} +  \sqrt{\frac{2\sigma^2}{\pi a} }\frac{ \exp( -\mu^2/2a\sigma^2)}{  \mbox{erfc}(-\mu/\sqrt{2a\sigma^2})},
   \end{equation}     
see \cite[eqn (26)]{DNR}. (The formula shown in \cite[eqn (21)]{DNR} has a misplaced $\sigma^2$ in its denominator; the formula in \cite[eqn (26)]{DNR} is correct.) 
 
\begin{lemma} \label{lem nu bounds}  Fix $\mu > 0$ and $a >0$ and $\sigma_0 > 0$.  Then the integrals 
   $\int_{\R_+} r^2 d\nu_{(\mu,a,\sigma)}(r)$ and $\int_{\R_+} r^{-2/3} d\nu_{(\mu,a,\sigma)}(r)$ are bounded above and the integral $\int_{\R_+} r^{2/3} d\nu_{(\mu,a,\sigma)}(r)$
is bounded away from $0$ for $0 < \sigma \le \sigma_0$. 
    \end{lemma}

\n{\bf Proof.}  For $\mu > 0$ we have $\mbox{erfc}(-\mu/\sqrt{2a\sigma^2})  >1$.  Using \eqref{int r2}we have 
   $$
   \int_{\R_+} r^2 d\nu_{(\mu,a,\sigma)}(r) \le \frac{\mu}{a} +  \sqrt{\frac{2\sigma^2}{\pi a} } \le \frac{\mu}{a} +  \sqrt{\frac{2\sigma_0^2}{\pi a} }$$
for all $\sigma \le \sigma_0$.  Next, if $r \le \sqrt{\mu/2a}$ then 
    $$\rho(r) \le \frac{2r\exp\left(-\frac{a}{2\sigma^2}\left[r^2-\frac{\mu}{a}\right]^2\right)}{\sqrt{\frac{\pi \sigma^2}{2a}}} \le \frac{2r}{\sqrt{\frac{\pi}{2a}}}\cdot \frac{e^{-\mu^2/8a\sigma^2}}{ \sigma} \le kr
    $$ 
for all $\sigma > 0$, where $k$ is a constant depending on $\mu$ and $a$.  Therefore for all $\sigma >0$
   $$\int_{\R_+} r^{-2/3} d\nu_{\mu,a,\sigma}(r) \le \int_0^{\sqrt{\mu/2a}} r^{-2/3} d\nu_{\mu,a,\sigma}(r) + \left(\frac{2a}{\mu}\right)^{1/3} 
    \le k \int_0^{\sqrt{\mu/2a}} r^{1/3}dr + \left(\frac{2a}{\mu}\right)^{1/3}.
   $$
Finally since $\mbox{erfc}(-\mu/\sqrt{2a\sigma^2})  < 2$ we have 
 \begin{align*}
  \int_{\R_+} r^{2/3}d\nu_{(\mu,a,\sigma)}(r)
  & \ge  \int_{\sqrt{\mu/a}}^\infty  r^{2/3} \cdot 
     \frac{2r\exp\left(-\frac{a}{2\sigma^2}\left[r^2-\frac{\mu}{a}\right]^2\right)}{2\sqrt{\frac{\pi \sigma^2}{2a}}}dr \\
     & =  \int_{\mu/a}^\infty
     s^{1/3} \cdot \frac{\exp\left(-\frac{a}{2\sigma^2}\left[s-\frac{\mu}{a}\right]^2\right)}{\sqrt{2 \pi \sigma^2/a}}ds\\
     & \ge \left(\frac{\mu}{a}\right)^{1/3} \int_{\mu/a}^\infty
     \frac{\exp\left(-\frac{a}{2\sigma^2}\left[s-\frac{\mu}{a}\right]^2\right)}{\sqrt{2 \pi \sigma^2/a}}ds
     = \frac{1}{2} \left(\frac{\mu}{a}\right)^{1/3},
     \end{align*}
and the proof is complete.  
        \qed

\begin{lemma} \label{lem nu bounds2}  Suppose $\mu > 0$.  For all measurable functions $f: \R \to [0,\infty)$ we have 
     \begin{equation} \label{normal}
     \int_{\R_+} f(r^2) d\nu_{(\mu,a,\sigma)}(r) \le 2 \E[f(\chi)]
     \end{equation}
where $\chi$ is a normal random variable with mean $\mu/a$ and variance $\sigma^2/a$.  In particular
    \begin{equation} \label{moment}
    \int_{\R_+} \left(r^2 - \frac{\mu}{a}\right)^2 d\nu_{(\mu,a,\sigma)}(r)\le \frac{2\sigma^2}{a}
    \end{equation}
and 
\begin{equation} \label{exp moment}
    \int_{\R_+} e^{\gamma(r^2-\mu/a)^2} d\nu_{(\mu,a,\sigma)}(r) \le \frac{2}{\sqrt{1-2\gamma \sigma^2/a}}
    \end{equation}   
for all $\gamma < a/(2\sigma^2)$.

\end{lemma} 
  
\n{\bf Proof.}  The inequality \eqref{normal} is proved by Chemnitz and Engel \cite[Lemma 3.11]{CE23}, and then \eqref{moment} and \eqref{exp moment} use standard facts about the normal random variable.  \qed

\section{Estimates involving $\int_M {\cal L}F(x)dm(x)$} \label{sec LF}

Suppose ${\cal L}$ is the generator of a diffusion process on a manifold $M$ with invariant probability measure $m$ and let $F \in C^2(M)$.  It is not always true that $\int_M {\cal L}F(x)dm(x) = 0$.  For example let ${\cal L}_r$ be the generator for $\{r_t: t \ge 0\}$ given by \eqref{r} above and let $F(r)= \log r$.  Then ${\cal L}F(r) = \mu - ar^2$, but $\int r^2 d\nu(r) \neq \mu/a$, see \eqref{int r2} above.  In order to apply the adjoint method on a non-compact space we will use the following strengthening of \cite[Prop 3]{BG01}.  

\begin{proposition} \label{prop-bg}
  Let $\{X_t:\ t\geq 0\}$ be a non-explosive 
diffusion process on a $\sigma$-compact manifold $M$ with
invariant probability measure $m$. Let ${\cal L}$ be an operator acting
on $C^2(M)$ functions that agrees with the generator of $\{X_t:\
t\geq 0\}$ on $C^2$ functions with compact support.  Let $F \in
C^2(M)$.  Suppose $F$ and ${\cal L}F$ are $m$ integrable.  Suppose also there exists a positive $V \in C^2(M)$ satisfying
${\cal L} V(x) \leq kV(x)$ for some $k<\infty$ such that $|F(x)|/V(x) \to 0$ as $|F(x)| \to \infty$. Then
 $$
 \int_M {\cal L}F(x)\, dm(x) = 0.
 $$
 \end{proposition}

 \n{\bf Proof.}  Let $\PP^x$ denote the law of the process
$\{X_t: t \ge 0\}$ with $X_0 = x$.  It\^o's formula implies that
 $
 N_t: = F(X_t) - F(x) - \int_0^t {\cal L}F(X_s)\,ds
 $
is a local martingale with respect to the measure $\PP^x$ for all $x
\in M$.  More precisely, if $\{K_n\}_{n=1}^\infty$ is a sequence
of compact subsets of $M$ such that $K_n\nearrow M$ and
$\delta_n=\inf \left\{t\geq 0:\ X_t \notin K_n\right\}$, then
$\delta_n \nearrow \infty$ and $N_{t \wedge \delta_n}$ is a
martingale with respect to the measure $\PP^x$ for all $x \in M$ and
all $n$.  In particular
  \begin{equation}
  \E^x(F(X_{t\wedge \delta_n})) - F(x) =
 \E^x\left( \int_0^{t \wedge \delta_n} {\cal L}F(X_s)ds \right). 
\label{mg}
 \end{equation}
 A similar calculation applied to the function $e^{-kt}V(y)$ shows
that $e^{-k(t\wedge \delta_n)} V(X_{t \wedge \delta_n})$ is a
$\PP^x$ supermartingale for all $x \in M$ and all $n$. It follows
that
 $$
  \E^x(V(X_{t \wedge \delta_n})) \le e^{kt} V(x). 
 $$
Given $\e > 0$ choose $c$ so that $|F(y)| \ge c$ implies $|F(y)|/V(y) <\e$.  Then 
    $$
    \E^x\left( |F(X_{t\wedge \delta_n})| 1_{|F(X_{t\wedge \delta_n})| \ge c}\right) \le \e \E^x(V(X_{t\wedge \delta_n})) \le \e e^{kt} V(x).
    $$  
Thus for each $t$ and $x$ the family $\{F(X_{t \wedge
\delta_n}): n \ge 1\}$ is $\PP^x$ uniformly integrable, and so 
    $$
     \E^x(F(X_{t\wedge \delta_n})) \to  \E^x(\lim_{n \to \infty}F(X_{t\wedge \delta_n})) = \E^x(F(X_t)).
     $$
Also
  $$
  \left|\int_0^{t \wedge \delta_n} {\cal L}F(X_s)ds \right| \le \int_0^t|{\cal L}F(X_s)|ds
  $$
and, using Fubini's theorem and the invariance of $m$,
  $$
 \int_M\left(\E^x \int_0^t |{\cal L}F(X_s)|ds \right)dm(x)
   = \int_0^t\left(\int_M \bigl(\E^x|{\cal L}F(X_s)|\bigr)dm(x)\right)ds =   t \int_M |{\cal L}F(x)|dm(x) <\infty 
   $$
for each $t > 0$. For fixed $t$ there exists a measurable subset $M_0  \subset M$ with $m(M_0) = 1$ such that     $$
   \E^x  \int_0^t |{\cal L}F(X_s)|ds <\infty
   $$
for all $x \in M_0$. Letting $n \to \infty$ in \eqref{mg} and using uniform integrability on the left and the dominated convergence theorem on the right, we get
 $$
 E^x(F(X_t)) - F(x) =
   E^x \left( \int_0^t {\cal L}F(X_s) ds \right) \quad \mbox{ for all }x \in M_0.
 $$
Integrating with respect to $m$ and using the invariance of $m$ again we obtain
  $$
 0 = \int_M  F(x) \,dm(x) -  \int_M F(x) \,dm(x) = t \int_M {\cal L}F(x)\,dm(x)
  $$
and the proof is complete. \hfill $\Box$

\bibliographystyle{plain}
\bibliography{shear.ref}

\begin{thebibliography}{10}

\bibitem{AX90}
S.~T. Ariaratnam and Wei~Chau Xie.
\newblock Lyapunov exponent and rotation number of a two-dimensional nilpotent
  stochastic system.
\newblock {\em Dynam. Stability Systems}, 5(1):1--9, 1990.

\bibitem{AIN04}
L.~Arnold, P.~Imkeller, and N.~Sri Namachchivaya.
\newblock The asymptotic stability of a noisy non-linear oscillator.
\newblock {\em J. Sound Vibration}, 269(3-5):1003--1029, 2004.

\bibitem{AOP86}
L.~Arnold, E.~Oeljeklaus, and \'{E}. Pardoux.
\newblock Almost sure and moment stability for linear {I}t\^{o} equations.
\newblock In {\em Lyapunov exponents ({B}remen, 1984)}, volume 1186 of {\em
  Lecture Notes in Math.}, pages 129--159. Springer, Berlin, 1986.

\bibitem{Arn98}
Ludwig Arnold.
\newblock {\em Random dynamical systems}.
\newblock Springer Monographs in Mathematics. Springer-Verlag, Berlin, 1998.

\bibitem{APW86}
Ludwig Arnold, George Papanicolaou, and Volker Wihstutz.
\newblock Asymptotic analysis of the {L}yapunov exponent and rotation number of
  the random oscillator and applications.
\newblock {\em SIAM J. Appl. Math.}, 46(3):427--450, 1986.

\bibitem{AM82}
E.~I. Auslender and G.~N. Mil'shtein.
\newblock Asymptotic expansions of the {L}iapunov index for linear stochastic
  systems with small noise.
\newblock {\em J. Appl. Math. Mech.}, 46:277--283, 1982.

\bibitem{Bax94}
Peter~H. Baxendale.
\newblock A stochastic {H}opf bifurcation.
\newblock {\em Probab. Theory Related Fields}, 99(4):581--616, 1994.

\bibitem{Bax97}
Peter~H. Baxendale.
\newblock Stability along trajectories at a stochastic bifurcation point.
\newblock In {\em Stochastic dynamics ({B}remen, 1997)}, pages 1--25. Springer,
  New York, 1999.

\bibitem{BG01}
Peter~H. Baxendale and Levon Goukasian.
\newblock Lyapunov exponents of nilpotent {I}t\^{o} systems with random
  coefficients.
\newblock {\em Stochastic Process. Appl.}, 95(2):219--233, 2001.

\bibitem{BG02}
Peter~H. Baxendale and Levon Goukasian.
\newblock Lyapunov exponents for small random perturbations of {H}amiltonian
  systems.
\newblock {\em Ann. Probab.}, 30(1):101--134, 2002.

\bibitem{BE23}
Maxime Breden and Maximilian Engel.
\newblock Computer-assisted proof of shear-induced chaos in stochastically
  perturbed {H}opf systems.
\newblock {\em Ann. Appl. Probab.}, 33(2):1052--1094, 2023.

\bibitem{CE23}
Dennis Chemnitz and Maximilian Engel.
\newblock Positive {L}yapunov exponent in the {H}opf normal form with additive
  noise.
\newblock {\em Comm. Math. Phys.}, 402(2):1807--1843, 2023.

\bibitem{CF98}
Hans Crauel and Franco Flandoli.
\newblock Additive noise destroys a pitchfork bifurcation.
\newblock {\em J. Dynam. Differential Equations}, 10(2):259--274, 1998.

\bibitem{DNR}
R.~E.~Lee DeVille, N.~Sri~Namachchivaya, and Zoi Rapti.
\newblock Stability of a stochastic two-dimensional non-{H}amiltonian system.
\newblock {\em SIAM J. Appl. Math.}, 71(4):1458--1475, 2011.

\bibitem{DFH08}
H.~A. Dijkstra, L.~M. Frankcombe, and A.~S. von~der Heydt.
\newblock A stochastic dynamical systems view of the {A}tlantic multidecadal
  oscillation.
\newblock {\em Philos. Trans. R. Soc. Lond. Ser. A Math. Phys. Eng. Sci.},
  366(1875):2545--2560, 2008.

\bibitem{DELR}
Thai~Son Doan, Maximilian Engel, Jeroen S.~W. Lamb, and Martin Rasmussen.
\newblock Hopf bifurcation with additive noise.
\newblock {\em Nonlinearity}, 31:4567--4601, 2018.

\bibitem{ELR}
Maximilian Engel, Jeroen S.~W. Lamb, and Martin Rasmussen.
\newblock Bifurcation analysis of a stochastically driven limit cycle.
\newblock {\em Comm. Math. Phys.}, 365(3):935--942, 2019.

\bibitem{FGS}
Franco Flandoli, Benjamin Gess, and Michael Scheutzow.
\newblock Synchronization by noise.
\newblock {\em Probab. Theory Related Fields}, 168(3-4):511--556, 2017.

\bibitem{GH90}
John Guckenheimer and Philip Holmes.
\newblock {\em Nonlinear oscillations, dynamical systems, and bifurcations of
  vector fields}, volume~42 of {\em Applied Mathematical Sciences}.
\newblock Springer-Verlag, New York, 1990.
\newblock Revised and corrected reprint of the 1983 original.

\bibitem{IK77}
Nobuyuki Ikeda and Shinzo Watanabe.
\newblock A comparison theorem for solutions of stochastic differential
  equations and its applications.
\newblock {\em Osaka Math. J.}, 14(3):619--633, 1977.

\bibitem{IL99}
Peter Imkeller and Christian Lederer.
\newblock An explicit description of the {L}yapunov exponents of the noisy
  damped harmonic oscillator.
\newblock {\em Dynam. Stability Systems}, 14(4):385--405, 1999.

\bibitem{IL01}
Peter Imkeller and Christian Lederer.
\newblock Some formulas for {L}yapunov exponents and rotation numbers in two
  dimensions and the stability of the harmonic oscillator and the inverted
  pendulum.
\newblock {\em Dyn. Syst.}, 16(1):29--61, 2001.

\bibitem{LSBY09}
Kevin~K. Lin, Eric Shea-Brown, and Lai-Sang Young.
\newblock Reliability of coupled oscillators.
\newblock {\em J. Nonlinear Sci.}, 19(5):497--545, 2009.

\bibitem{LY08}
Kevin~K. Lin and Lai-Sang Young.
\newblock Shear-induced chaos.
\newblock {\em Nonlinearity}, 21(5):899--922, 2008.

\bibitem{LWY13}
Kening Lu, Qiudong Wang, and Lai-Sang Young.
\newblock Strange attractors for periodically forced parabolic equations.
\newblock {\em Mem. Amer. Math. Soc.}, 224(1054):vi+85, 2013.

\bibitem{MTIII}
Sean~P. Meyn and R.~L. Tweedie.
\newblock Stability of {M}arkovian processes. {III}. {F}oster-{L}yapunov
  criteria for continuous-time processes.
\newblock {\em Adv. in Appl. Probab.}, 25(3):518--548, 1993.

\bibitem{PRK}
Arkady Pikovsky, Michael Rosenblum, and J\"{u}rgen Kurths.
\newblock {\em Synchronization: A Universal Concept in Nonlinear Sciences}.
\newblock Cambridge Nonlinear Science Series. Cambridge University Press, 2001.

\bibitem{PW88}
M.~A. Pinsky and V.~Wihstutz.
\newblock Lyapunov exponents of nilpotent {I}t\^{o} systems.
\newblock {\em Stochastics}, 25(1):43--57, 1988.

\bibitem{CTDN-II}
A.~Tantet, M.~Chekroun, H.~Dijkstra, and D.~Neelin.
\newblock {R}uelle-{P}ollicott resonances of stochastic systems in reduced
  state space. {P}art {II}: {S}tochastic {H}opf bifurcation.
\newblock {\em J. Stat. Phys.}, 179:1403--1448, 2020.

\bibitem{WY03}
Qiudong Wang and Lai-Sang Young.
\newblock Strange attractors in periodically-kicked limit cycles and {H}opf
  bifurcations.
\newblock {\em Comm. Math. Phys.}, 240(3):509--529, 2003.

\bibitem{Wie}
Sebastian Wieczorek.
\newblock Stochastic bifurcation in noise-driven lasers and {H}opf oscillators.
\newblock {\em Phys. Rev. E (3)}, 79(3):036209, 10, 2009.

\end{thebibliography}

\end{document}